\documentclass{elsart}
\usepackage{natbib}

\usepackage{amsmath}
\usepackage{amssymb}
\newcommand{\Q}{{\mathbb{Q}}}
\newcommand{\R}{{\mathbb{R}}}
\newcommand{\Z}{{\mathbb{Z}}}
\newcommand{\N}{{\mathbb{N}}}
\newcommand{\C}{{\mathbb{C}}}
\newcommand{\A}{{\mathbb{A}}}
\newcommand{\Oh}{{\mathcal{O}}}

\newcommand{\arr}{{\,\Longrightarrow\,}}
\newcommand{\abs}[1]{\lvert#1\rvert}

\newtheorem{lmm}{Lemma}

\newtheorem{example}{Example}

\theoremstyle{remark}
\newtheorem{remark}{Remark}

\numberwithin{equation}{section}

\begin{document}

\begin{frontmatter}

\title{Relation ideals and the Buchberger--M\"oller algorithm}

\author{Mathias Lederer}
\ead{mlederer@mathematik.uni-bielefeld.de}

\address{Fakult\"at f\"{u}r Mathematik, Universit\"{a}t Bielefeld, Bielefeld, Germany}

\begin{abstract}
  We construct a Gr\"{o}bner Basis of the relation ideal of a polynomial, give an interpolation formula 
  for the basis elements and explain the connection of the interpolation formula to the Buchberger--M\"oller 
  algorithm. We present a situation in which the usage of the Buchberger--M\"oller algorithm is 
  obsolete since one can compute its result directly. We prove a constructive version of a theorem of Galois, 
  concerning the solvability of rational polynomials of prime degree. 
  Computations are carried out for a number of example polynomials. 
\end{abstract}


\end{frontmatter}

\section{Introduction}\label{secintro}

Let $K$ be a field and $f=Z^n+a_1Z^{n-1}+\ldots+a_n$ a monic, irreducible, and separable polynomial in the 
variable $Z$ over $K$. Let $x=(x_1,\ldots,x_n)$ be the $n$-tuple of the zeros of $f$ in some 
field extension of $K$, and let $T=(T_1,\ldots,T_n)$ be indeterminates over $K$. We consider the set
\begin{equation}\label{ideal}
  I=\{P \in K[T];P(x)=0\}\,,
\end{equation}
which is an ideal in $K[T]$, called the \textit{relation ideal of $f$}. 
We studied this ideal already in the paper \cite{parma}; nevertheless, we will repeat 
some of the arguments and constructions of \cite{parma} in the present paper, which may need some 
more explanation than was given in the proceedings paper \cite{parma}. 
In this section, we will start by providing the motivation why it is important to know the ideal $I$ well. 
Furthermore, we will mention some results prior to our own work. 

Let $L=K(x_1,\ldots,x_n)$ be the splitting field of $f$. Suppose we have a method at hand that allows us to 
do computations in the field $K$, such as trivially in the case $K=\Q$, or in the case that $K$ 
is a number field $K=\Q(\alpha)$ given by the minimal polynomial of $\alpha$. 
Given such a ground field $K$, it is not clear how to do computations in the splitting field $L$ of $f$. 
The only case in which this is clear is the special case where $L=K(x_i)$ 
for some zero $x_i$ of $f$. In this case we will make use of the isomorphism $K[Z]/(f)\xrightarrow{\sim} L$ 
that sends the residue class of $P$ to $P(x)$ -- we know how to do computations in the quotient $K[Z]/(f)$, 
so by virtue of the isomorphism, we also know how to do computations in $L$. 
In the general case, that is the case when $L=K(x_1,\ldots,x_n)$ is strictly bigger than some $K(x_i)$, 
there exists an analogous isomorphism:
\begin{thm}\label{proiso}
  The map $\phi:K[T]/I \to L: P\mapsto P(x)$ is a $K$-algebra isomorphism.
\end{thm}
This is an easy consequence of the Homomorphism Theorem. The isomorphism of the theorem allows 
us to do computations in the field $L=K(x_1,\ldots,x_n)$ -- but only if we know how to do 
computations in the quotient $K[T]/I$! Since we have assumed complete knowledge about computing in $K$, 
knowledge about computing in $K[T]/I$ is equivalent to knowledge of a Gr\"{o}bner basis of $I$. 

We can rephrase this as follows: The question, ``How can one do computations in the splitting field of a polynomial?'' 
is at least as old as field theory itself. If we had a Gr\"{o}bner basis of the relation ideal $I$, 
this question would be solved (providing that we know how to do computations in the ground field). 
However, the definition of $I$ does not automatically lead us to a generating set of $I$. 
Several people have taken different approaches to this problem, which we will shortly describe. 

One approach is \cite{anai}. The subject of this paper is the computation of the Galois group of a polynomial $f$. 
The authors distinguish between two ways of computing the Galois group, \textit{table-based methods} 
on the one hand and \textit{direct methods} on the other. We shall not explain 
the differences between these two sets of methods but merely explain the information they give us. 
The result of a table-based method is a conjugacy class of subgroups of the symmetric group $S_n$ 
(remember that $n$ is the degree of $f$), 
having the property that, given the zeros of $f$, there is a subgroup $G$ of $S_n$ 
in the given conjugacy class such that $G$ can be identified with the Galois group of $f$ via 
$\sigma(x_i)=x_{\sigma(i)}$ for all $\sigma\in G$ and all $i\in\{1,\ldots,n\}$. 
A direct method provides us with more information: It computes not only a conjugacy class 
of subgroups of $S_n$ but one particular $G$, as above, for one given numbering $x=(x_1,\ldots,x_n)$ 
of the zeros. When we speak of zeros in this context, we mean approximations of the numerical
values of the zeros. Think of the case $K=\Q$, then $x=(x_1,\ldots,x_n)$ is an $n$-tuple of complex or
$p$-adic numbers approximating the complex or $p$-adic zeros of $f$. 

Let us recall where the ambiguity in the description of the Galois group comes from. 
Suppose we are given a subgroup $G$ of $S_n$ such that the Galois group of $f$ can be identified 
with $G$ via $\sigma(x_i)=x_{\sigma(i)}$ for all $\sigma\in G$ and all $i\in\{1,\ldots,n\}$. 
For this description, the knowledge of the numbering $x=(x_1,\ldots,x_n)$ of the zeros of $f$ is 
necessary. If we chose a different numbering of the zeros, thus if we replaced $x=(x_1,\ldots,x_n)$ by 
$y=(y_1,\ldots,y_n)$, where $y_i:=x_{\tau(i)}$ for some $\tau\in S_n$, the subgroup $H$ of $S_n$ such that 
$\sigma(y_i)=y_{\sigma(i)}$ for all $\sigma\in H$ would be $H=\tau^{-1}G\tau$. From this follows in particular that
if we do not fix one particular numbering of the zeros but only look at the set of zeros, 
we can never know the Galois group of $f$ any better than up to conjugacy in $S_n$. 

The subject of \cite{anai} is a direct method for the computation of the Galois group of a polynomial. 
A Gr\"obner basis of $I$ is constructed as a by-product. The idea of this method is 
to successively factor $f$ over certain finite extensions of $K$. 
The elements of the Gr\"{o}bner basis of $I$ are derived from the factorisations of $f$. 
However, factorising $f$ over an algebraic extension of $K$, as used in the method of \cite{anai}, 
is not an easy task. Thus, one would like to obtain the generators of the relation ideal in a 
more direct way than indicated in \cite{anai}.

Paper \cite{mckay} takes a step towards a direct computation of the generators of the relation ideal:
Certain polynomials $\Gamma_i$, $i=1,\ldots,n$, with coefficients in $K$ are constructed. 
The ingredients for the computation of $\Gamma_i$ are the same as the ingredients for our computation 
in \cite{parma} and in the present paper: We take the Galois group $G$ of $f$ to be given, 
more precisely as a subgroup of $S_n$ permuting the zeros $x=(x_1,\ldots,x_n)$ of $f$, 
along with approximations of the zeros of $f$. (Thus we take as given the data which a direct method
for the computation of the Galois group gives us, in the terminology of \cite{anai}.)
The virtue of the polynomials $\Gamma_i$ is as follows: It may and will happen that the splitting field of $f$ 
is obtained by adjoining \textit{less than} $n$ zeros to the ground field $K$, 
say, $L=K(x_1,\ldots,x_k)$ for $k<n$ (in fact, $k\leq n-1$ always holds).
In particular, it is possible to write the zeros $x_i$, $i>k$, as a rational expression in $x_1,\ldots,x_k$
with coefficients in $K$. This expression can be gleaned from the polynomials $\Gamma_i$, $i>k$. 
Thus by \cite{mckay}, some of the relations between the roots of $f$ are known. 
The polynomials $\Gamma_i$ are given by an explicit formula involving $G$ and $x$ and are in practice 
(over $K=\Q$) computed by using complex approximations of $x$. In the last section of \cite{mckay}, 
the authors state that it is natural to ask how one might find all the rational relations between the roots of $f$. 

This question was answered in our paper \cite{parma} and is now explained at more detail. 
To give a rough sketch, we will proceed as follows: We will construct a Gr\"obner basis of the relation ideal $I$. 
The basic idea will be to use the same finite extensions of $K$ as were used in \cite{anai}. 
We will define them in Section \ref{secgen} below. However, we will not obtain the generators 
by factorising $f$ over the extensions of $K$, as was done in \cite{anai}. 
Instead, we will suppose that we already have the Galois group $G$ of $f$ 
as a subgroup of $S_n$ permuting the (approximations of the) zeros of $f$. 
From that we will compute the elements of the Gr\"{o}bner basis of $I$ directly, 
giving an interpolation formula in the spirit of Lagrange interpolation. 
This is the formula \eqref{lag}, see Theorem \ref{thmlag} in Section \ref{secint}. 
Our interpolation formula will be structurally similar to the formulas for the polynomials $\Gamma_i$ of \cite{mckay}. 
In particular, the formulas of \cite{mckay} and our interpolation formula will use the same ingredients, 
i.e. the zeros of $f$ and the Galois group $G$, permuting the zeros. 

Only when we have proved the interpolation formula for the generators, 
a sufficient degree of explicitness is reached so that we can consider the question,
``What are the generators of the relation ideal?'' answered. Thus also the question,
``How can one do computations in the splitting field of a polynomial?'' is solved. 

In Section \ref{secmoe} we will explain that our interpolation formula does not exist independently of 
other ideas in mathematics. In fact our interpolations formula is exactly the formula that comes out when 
we apply the Buchberger--M\"oller algorithm (see \cite{buchbergermoeller}) 
to the interpolation problem of Theorem \ref{procoe}. Of course, one need not evoke the entire
Buchberger--M\"oller algorithm in order to find the generators of the relation ideal, 
since we already know the interpolation formula \eqref{lag} for the generators. 
So it is only natural to ask for an interpolation formula like ours for a more general interpolation problem. 
In Section \ref{secmoe} we also take this question a bit further: 
We observe that our interpolation formula merely uses the transitive and faithful operation of the Galois group, 
that permutes the zeros of $f$. Hence an analogous formula like ours holds true also when the group with 
these properties is not necessarily a Galois group that permutes the zeros of a polynomial. 

In Section \ref{secgal} we apply our results to a classical theorem of Galois, as to be found in \cite{cohn} 
or \cite{hup}:

\begin{quotation}
  A rational polynomial $f$ of degree $p$, where $p$ is a prime number, is solvable by radicals if and only if each 
  zero of $f$ can be expressed as a polynomial (with rational coefficients) in any two other zeros.
\end{quotation}

\noindent As Olaf Neumann commuticated personally (to Kurt Girstmair), this theorem was crucial for Galois theory 
to be accepted in the mathematical community, since its statement was easy to understand for any mathematician 
around 1840, yet its proof requires Galois theory, still a new tool at that time. The theorem asserts the existence 
of a rational polynomial without giving its precise form. Also, this polynomial was never constructed 
in the literature. We will see that this polynomial is one of the generators of the relation ideal that we had 
constructed before. In this sense our paper gives a constructive version of an important existence theorem. 

In Section \ref{secnum} we make use of some $p$-adic techniques (similar to those of \cite{klue}) 
to compute a number of examples over the ground field $\Q$ in Section \ref{secexa}. 
However, using complex approximations of the zeros might just as well serve to treat the examples. 
In Section \ref{secque}, we point at a property of the generators that cannot be explained 
within the theory that was used in this paper. Furthermore, we will try to entice the reader
to study our concept of interpolation and the Buchberger--M\"oller algorithm in more generality. 

\section{Generators of the relation ideal}\label{secgen}

For $i=1,\ldots,n$, define the field $K_i=K(x_1,\ldots,x_i)$, and set $K_0=K$. Then clearly $K_i=K_{i-1}(x_i)$
for $i=1,\ldots,n$, and $x_i$ is a primitive element of the field $K_i$ over the field $K_{i-1}$. 
Let $f_i$ be the minimal polynomial of $x_i$ over $K_{i-1}$. Then $d_i=\deg(f_i)$ is the degree
of the field extension $K_i|K_{i-1}$. Therefore the polynomial $f_i$ has the shape
\begin{equation*}
  f_i=T_i^{d_i}+\sum_{k_i=1}^{d_i}b_{i,k_i}T_i^{d_i-k_i}\,,
\end{equation*}
where all coefficients $b_{i,k_i}$ lie in $K_{i-1}$. Since the degree of the field extension 
$K_i|K_{i-1}$ equals $d_i$, the degree formula shows that the degree of the field extension $K_i|K$ equals 
$d_1\ldots d_i$. It is easy to see that the family $x_1^{d_1-k_1}\ldots x_{i-1}^{d_{i-1}-k_{i-1}}$, where 
$1\leq k_j\leq d_j$ for $j=1,\ldots,i-1$, is a $K$-basis of $K_{i-1}$. 
Thus the coefficients of the polynomial $f_i$ can uniquely be written as
\begin{equation*}
  b_{i,k_i}=\sum_{k_1=1}^{d_1}\ldots \sum_{k_{i-1}=1}^{d_{i-1}}b_{i,k_1,\ldots,k_i}
  x_1^{d_1-k_1}\ldots x_{i-1}^{d_{i-1}-k_{i-1}}\,,
\end{equation*}
all $b_{i,k_1,\ldots,k_i}$ belonging to $K$. (For $i=1$, no summation needs to be done and we simply have 
$b_{1,k_1}\in K$.) We obtain:
\begin{equation}\label{fplain}
  f_i=T_i^{n_i}+\sum_{k_1=1}^{d_1}\ldots\sum_{k_i=1}^{d_i}
  b_{i,k_1,\ldots,k_i}x_1^{d_1-k_1}\ldots x_{i-1}^{d_{i-1}-k_{i-1}}T_i^{d_i-k_i}\,.
\end{equation}
Now we define
\begin{equation}\label{fhat}
  \widehat{f_i}=T_i^{d_i}+\sum_{k_1=1}^{d_1}\ldots\sum_{k_i=1}^{d_i}
  b_{i,k_1,\ldots,k_i}T_1^{d_1-k_1}\ldots T_{i-1}^{d_{i-1}-k_{i-1}}T_i^{d_i-k_i}\,,
\end{equation}
thus $\widehat{f_i}\in K[T_1,\ldots,T_i]$. We will use the identities $f_i=\widehat{f_i}(x_1,\ldots,x_{i-1},T_i)$ and 
$f_i(x_i)=\widehat{f_i}(x_1,\ldots,x_{i-1},x_i)$ later. In what follows all polynomials 
$\widehat{f_1}, \ldots,\widehat{f_i}$ are considered to lie in $K[T_1,\ldots,T_i]$. 

\begin{thm}\label{thmgen}
  The $K$-algebras $K[T_1,\ldots,T_i]/(\widehat{f_1},\ldots,\widehat{f_i})$ and $K_i$ are 
  isomorphic via the map $\phi_i$ which assigns to $P$ the value $P(x_1,\ldots,x_i)$. 
\end{thm}

\begin{pf}
  Consider the evaluation $\psi_i:K[T_1,\ldots,T_i]\to K(x_1,\ldots,x_i)$ defined by 
  $\psi_i(P)=P(x_1,\ldots,x_i)$. 
  We have to show that $\ker(\psi_i)=(\widehat{f_1},\ldots,\widehat{f_i})$. Obviously, 
  $\ker(\psi_i)\supseteq(\widehat{f_1},\ldots,\widehat{f_i})$. We show the converse inclusion by induction over $i$. 

  For $i=1$, the assertion is well known: The mapping $K[T_1]\to K(x_1):P\mapsto P(x_1)$ is an isomorphism. 
  For $i>1$, we will make use of two isomorphisms, the first is an analogue of the one we just used, that is 
  the isomorphism
  \begin{equation*}
    \alpha:K_{i-1}[T_i]/(f_i)\to K_i:P\mapsto P(x_i)\,.
  \end{equation*}
  For the second, the induction hypothesis says that 
  \begin{equation*}
    \phi_{i-1}:K[T_1,\ldots,T_{i-1}]/(\widehat{f_1},\ldots,\widehat{f_{i-1}})\to K_{i-1}:
    P\mapsto P(x_1,\ldots,x_{i-1})
  \end{equation*}
  is an isomorphism. We adjoin to the domain of definition of $\phi_{i-1}$ and to the range of $\phi_{i-1}$ 
  the variable $T_i$ and obtain the second isomorphism,
  \begin{equation*}
    \beta:K[T_1,\ldots,T_i]/(\widehat{f_1},\ldots,\widehat{f_{i-1}})\to K_{i-1}:
    P\mapsto P(x_1,\ldots,x_{i-1},T_i)\,.
  \end{equation*}

  Let $P=P(T_1,\ldots,T_i)$ lie in the kernel of $\psi_i$. In view of $\alpha$, we conclude that 
  $P(x_1,\ldots,x_{i-1},T_i)$ is a multiple of $f_i$ by a polynomial in $K_{i-1}[T_i]$, so 
  \begin{equation*}
    P(x_1,\ldots,x_{i-1},T_i)=Q(x_1,\ldots,x_{i-1},T_i)f_i(T_i)\,,
  \end{equation*}
  for a suitable $Q\in K[T_1,\ldots,T_i]$. In view of $\beta$, the equation 
  \begin{equation*}
    P(x_1,\ldots,x_{i-1},T_i)-Q(x_1,\ldots,x_{i-1},T_i)f_i(T_i)=0\,,
  \end{equation*}
  shows that $P-Q\widehat{f_i}$ lies in the ideal spanned by $\widehat{f_1},\ldots,\widehat{f_{i-1}}$. 
  This shows that $P$ lies in the ideal spanned by $\widehat{f_1},\ldots,\widehat{f_i}$. 
\end{pf}

\begin{remark}
  It is quite natural to build up the splitting field $L$ of $f$ by the chain 
  \begin{equation}\label{chain}
    K_0\subset K_1\subset\ldots\subset K_n=L
  \end{equation}
  of monogenic field extensions, as we have done above. In \cite{gro} and later in \cite{anai}, 
  the same intermediate fields of $K$ and $L$ are used, and $K_i$ is presented as in Theorem \ref{thmgen}. 
  In particular, the existence of the Gr\"obner basis of the relation ideal of Theorem \ref{thmgen} was 
  -- at least in principle -- known to Wolfgang Gr\"obner himself. 
  However, his approach is not exactly the same the approach chosen here, so I have decided to present this 
  (shorter) proof of Theorem \ref{thmgen} instead of just giving a reference. 

  Gr\"obner (see \cite{gro}) defines the generators $\widehat{f_i}$ as we did, 
  that is by requiring $f_i$ to be the minimal polynomial of $x_i$ over $K_{i-1}$, 
  whereas Anai, Noro and Yokoyama (see \cite{anai}) require $f_i$ to be 
  a certain irreducible factor over $K_{i-1}$ of $f$. Of course these two characterisations are equivalent. 
  In \cite{anai}, the authors use their characterisation of $f_i$ in order to effectively compute 
  these polynomials. So they split $f$ into irreducible factors over certain field extensions of $K$. 
  In practice, this is not an easy task and will slow down the algorithm. 
  We will overcome this obstacle by using the interpolation 
  formula for $\widehat{f_i}$, which we will prove in the next section. 
\end{remark}

\begin{remark}
  Aubry and Valibouze extended Gr\"obner's work to a more general situation in \cite{aub}. 
  They proved the existence of a basis having the property of being ``separable triangular'' 
  (see \cite{aub} for the definition) for a class of ideals to which the relation ideal belongs. 
  Theorem \ref{thmgen} is basically a special case of their main theorem. 

  Let us adopt some of their notations, for the time being. A set $H$ of polynomials in $K[T_1,\ldots,T_i]$ is called
  \textit{triangular} if $H=\{H_1(T_1),\ldots,H_i(T_1,\ldots,T_i)\}$, where
  the $j$-th polynomial $H_j$ is monic as a polynomial in $T_j$ with $\deg_j(H_j)>0$. Here, and in the rest of 
  this paper, $\deg_j$ denotes the degree of $h$ in $T_j$. As the authors of \cite{aub} remark, 
  if $H$ is a triangular set then the ideal spanned by $H$ in $K[T_1,\ldots,T_i]$ is zero-dimensional 
  and $H$ is a reduced Gr\"{o}bner basis of the corresponding ideal, where we use the lexicographic 
  ordering with $T_1<\ldots<T_n$ on $K[T]$. (See \cite{beck}, \cite{buch}, or \cite{gro}.) 
  The property of being triangular is particularly valuable for practical purposes since the reduction of a 
  polynomial $P\in K[T]$ modulo $H$ comes down to dividing $P$ successively by the monic polynomials 
  $H_1,\ldots,H_i$ (in any order). By construction the set $\{\widehat{f_1},\ldots,\widehat{f_i}\}$ is clearly 
  a triangular set in $K[T_1,\ldots,T_i]$, thus a Gr\"obner basis of the corresponding ideal. 
\end{remark}

\section{The vanishing property characterising the generators}\label{secvan}

Although the generators $\widehat{f_i}$ are unambiguously defined and we even know that they have the shape 
\eqref{fhat} with all coefficients $b_{i,k_1,\ldots,k_i}$ lying in $K$, we do not know the numerical values of 
these coefficients. The reason is that, if given only the polynomial $f$, we do not have all the 
minimal polynomials $f_i$ at hand. In this section we will prove that every $\widehat{f_i}$ vanishes on 
a certain subset of $L^i$, and that this vanishing property, together with the constraint \eqref{fhat} 
on the degree of $\widehat{f_i}$, characterises $\widehat{f_i}$. 

In this section and onwards, we will make use of the Galois group $G$ of the field extension $L|K$. 
We assume that we have $G$ given as a subgroup of $S_n$ such that $\sigma(x_i)=x_{\sigma(i)}$ for all $i$. 
Furthermore, we will make use or the subgroup $G_i$ of $G$, which we define to be the group corresponding to 
the intermediate field $K_i$ of $L|K$. By definition of $K_i$, we have
$G_i=\{\sigma\in G;\sigma(x_j)=x_j,\text{ for all }j\leq i\}$. 

\begin{lmm}\label{lmmkli}
  Let $L|K$ and $G$ be as above. For $y\in L$, define 
  \begin{equation*}
    Gy= \begin{pmatrix} \sigma_1(y) \\ \vdots \\ \sigma_N(y) \end{pmatrix}\in L^N\,.
  \end{equation*}
  Then for arbitrary $y_1,\ldots,y_r\in L$, the following statements are equivalent:
  \begin{itemize}
  \item[(i)] $y_1,\ldots,y_r\in L$ are $K$-linearly independent.
  \item[(ii)] $Gy_1,\ldots,Gy_r\in L^N$ are $L$-linearly independent.
  \end{itemize}
\end{lmm}

\begin{pf}
  We point out that this lemma is quite similar to Artin's Lemma and only give a proof of the nontrivial direction 
  $(i)\arr(ii)$. Assume that $y_1,\ldots,y_r$ are $K$-linearly independent but $Gy_1,\ldots,Gy_k$ (where $k<r$) 
  form an $L$-basis of $_L\langle Gy_1,\ldots,Gy_r\rangle$. Then there exist uniquely determined coefficients 
  $\lambda_1,\ldots,\lambda_r\in L$ satisfying $Gy_{k+1}=\sum_{i=1}^{k}\lambda_iGy_i$. For every $\sigma\in G$, 
  there exists a matrix $P\in GL_N(K)$ satisfying $\sigma(Gy)=PGy$ for all $y\in L$. We obtain 
  $PGy_{k+1}=\sigma(y_{k+1})=\sum_{i=1}^{k}\sigma(\lambda_i)\sigma(Gy_i)
  =\sum_{i=1}^{k}\sigma(\lambda_i)PGy_i=P\sum_{i=1}^{k}\sigma(\lambda_i)Gy_i$ and therefrom $Gy_{k+1}=
  \sum_{i=1}^{k}\sigma(\lambda_i)Gy_i$. Since $Gy_1,\ldots,Gy_k$ is a basis of 
  $_L\langle Gy_1,\ldots,Gy_r\rangle$, $Gy_{k+1}$ is uniquely written as an $L$-linear combination of these vectors, 
  and therefore $\sigma(\lambda_i)=\lambda_i$ for $i=1,\ldots,k$ and for all $\sigma\in G$. Since $K$ is the 
  fixed field of $G$, the coefficient $\lambda_i$ must lie in $K$ for all $i=1,\ldots,k$. Thus also $y_{k+1}$ lies in 
  $_K\langle y_1,\ldots,y_r\rangle$, a contradiction to the $K$-linear independence of $y_1,\ldots,y_r$.
\end{pf}

\begin{thm}\label{procoe}
  $L[T_1,\ldots,T_i]$ contains exactly one polynomial of the shape \eqref{fhat} vanishing at 
  $(\sigma(x_1),\ldots,\sigma(x_i))$ for all $\sigma\in G$. The coefficients of this polynomial lie in $K$.
\end{thm}

\begin{pf}
  First we note that $\widehat{f_i}$ has the desired property: $f_i(x_i)=\widehat{f_i}(x_1,\ldots,x_i)=0$, 
  and also $\sigma(\widehat{f_i}(x_1,\ldots,x_i))=\widehat{f_i}(\sigma(x_1),\ldots,\sigma(x_i))=0$ for all 
  $\sigma\in G$. This proves the existence as claimed in the Theorem. Of course, the coefficients of 
  $\widehat{f_i}$ lie in $K$. 

  Since the family $x_1^{d_1-k_1}\ldots x_i^{d_i-k_i}$, where $1\leq k_j\leq d_j$ for $j=1,\ldots,i$, 
  is a $K$-basis of $K_i$, Lemma \ref{lmmkli} implies that the family $(Gx_1^{d_1-k_1}\ldots x_i^{d_i-k_i})$, 
  where $1\leq k_j\leq d_j$ for $j=1,\ldots,i$, is $L$-linearly independent. Thus the coefficients 
  $b_{i,k_1,\ldots,k_i}\in L$ in the sum 
  \begin{equation*}
    -Gx_i^{d_i}=\sum_{k_1=1}^{d_1}\ldots \sum_{k_i=1}^{d_i}b_{i,k_1,\ldots,k_i}Gx_1^{d_1-k_1}\ldots x_i^{d_i-k_i}
  \end{equation*}
  are uniquely determined. In other words, the coefficients of a polynomial having the shape \eqref{fhat} 
  are uniquely determined under the assumption that the polynomial vanishes at $(\sigma(x_1),\ldots,\sigma(x_i))$ 
  for all $\sigma\in G$. This proves the uniqueness as claimed in the theorem. 
\end{pf}

\section{An interpolation formula for the generators}\label{secint}

After Theorem \ref{procoe} there remains the task of finding a polynomial with coefficients in $L$ vanishing in all 
$(\sigma(x_1),\ldots,\sigma(x_i))$ and having the degree as described in \eqref{fhat}. 
This polynomial will be $\widehat{f_i}$. In this section we present an interpolation formula for $\widehat{f_i}$, 
which strongly reminds us of Lagrange interpolation. We will discuss the relation of our interpolation formula
to the results of \cite{anai} and \cite{mckay}. 

For the interpolation we will need the following sets: For $\rho\in G$ and $i=1,\ldots,n$, define
\begin{equation*}
  B{(\rho,i)}=\{\sigma(x_i);\sigma\in G,\sigma|_{K_{i-1}}=\rho|_{K_{i-1}}\}\setminus\{\rho(x_i)\}
\end{equation*}
In other words: $B{(\rho,i)}$ consists of the translates $\sigma(x_i)$, where $\sigma$ runs through all 
extensions of $\rho|_{K_{i-1}}$ to $K_i$, minus the element $\rho(x_i)$. 

\begin{lmm}\label{lmmbri}
  $|B{(\rho,i)}|=d_i-1$ for all $\rho\in G$ and for all $i\in\{1,\ldots,n\}$.
\end{lmm}

\begin{pf}
  The number of extensions $\sigma$ of $\rho|_{K_{i-1}}$ to $K_i$ equals the degree of the field extension 
  $K_i|K_{i-1}$, i.e. $d_i$. Two extensions of this kind are different if and only if they take different values 
  $\sigma(x_i)$, since $x_i$ generates $K_i$ over $K_{i-1}$. 
\end{pf}

\begin{thm}\label{thmlag}
  The $i$-th generating polynomial $\widehat{f_i}$ of the relation ideal $I$ is given by
  \begin{equation}\label{lag}
    \widehat{f_i}=T_i^{d_i}-\sum_{\rho\in G/G_i}\rho(x_i)^{d_i}
    \prod_{y_1\in B{(\rho,1)}}\frac{T_1-y_1}{\rho(x_1)-y_1}\ldots
    \prod_{y_i\in B{(\rho,i)}}\frac{T_i-y_i}{\rho(x_i)-y_i}\,.
  \end{equation}
\end{thm}

\begin{pf}
  We define $g$ by the right hand side of \eqref{lag} and prove $\widehat{f_i}=g$. Lemma \ref{lmmbri} 
  shows that $\text{deg}_j(g)\leq d_j-1$, for $j=1,\ldots,i-1$. Clearly $\text{deg}_i(g)=d_i$. 
  Thus the multidegree of $g$ has the properties that we demanded for $\widehat{f_i}$. 
  I we can prove that $g(\sigma(x_1),\ldots,\sigma(x_i))=0$ for all $\sigma\in G$, it will follow from 
  Theorem \ref{procoe} that $\widehat{f_i}$ and $g$ coincide. 

  So let $\sigma\in G$ be given. Take $\rho^{\prime}\in G/G_i$ such that $\sigma=\rho^{\prime}\tau$, for a suitable 
  $\tau\in G_i$. In particular, for $j=1,\ldots,i$ we have $\sigma(x_j)=\rho^{\prime}\tau(x_j)=\rho^{\prime}(x_j)$ 
  since $\tau(x_j)=x_j$. In order to show that $g(\sigma(x_1),\ldots,\sigma(x_i))=0$, we focus our attention 
  on the sum $\sum_{\rho\in G/G_i}$. The automorphisms $\rho$ occurring as summation index belong to the two 
  categories $\rho=\rho^{\prime}$ and $\rho\neq\rho^{\prime}$. If $\rho=\rho^{\prime}$, the respective summand of 
  becomes $\sigma(x_i)^{d_i}$ for in this case $\sigma(x_j)=\rho^{\prime}(x_j)=\rho(x_j)$, hence 
  $(\sigma(x_j)-y_j)/(\rho(x_j)-y_j)=1$ for all $j=1,\ldots,i$. In the case $\rho\neq\rho^{\prime}$ we can 
  find a number $j\in\{1,\ldots,i\}$ satisfying $\rho(x_j)\neq\rho^{\prime}(x_j)$. Thus $\rho^{\prime}(x_j)$ 
  lies in $B{(\rho,j)}$, and therefore there is a $y_j\in B{(\rho,j)}$ such that $y_j=\rho^{\prime}(x_j)=\sigma(x_j)$. 
  The summand corresponding to this $\rho$ vanishes, since the product occurring in in the sum 
  contains the factor $\sigma(x_j)-y_j$ where $y_j=\rho^{\prime}(x_j)=\sigma(x_1)$. Altogether, we obtain
  $g(\sigma(x_1),\ldots,\sigma(x_i))=\sigma(x_i)^{d_i}-\sigma(x_i)^{d_i}=0$. 
\end{pf}

\begin{remark}
  We discuss the relation of our work to \cite{mckay}. The central notion of \cite{mckay} is that of a 
  \textit{minimal strong base} for $f$ over $K$. This is is a set $\{x_{i_1},\ldots,x_{i_k}\}$ of roots 
  of $f$ such that $L=K(x_{i_1},\ldots,x_{i_k})$ and if $L=\Q(x_{j_1},\ldots,x_{j_m})$, then $k\leq m$. 
  Of course, one can always choose a numbering of the zeros of $f$ such that $x_1,\ldots,x_k$ is 
  a minimal strong base. 

  Given a minimal strong base, one can write all other zeros $x_i$, $i>k$, as a rational expression over $K$ 
  in the elements of the minimal strong base. Stauduhar and McKay (see \cite{mckay}) 
  use the following strategy in order to obtain this expression. They consider the polynomials 
  \begin{equation}\label{mckay}
    \Gamma_i=\sum_{\rho\in G}\rho(x_i)\frac{f(T_1)}{T_1-\rho(x_1)}\ldots\frac{f(T_k)}{T_k-\rho(x_k)}\,,
  \end{equation}
  which have coefficients in $K$, and satisfy $\Gamma_i(x_1,\ldots,x_k)=x_if^\prime(x_1)\ldots f^\prime(x_k)$. 
  Thus one can compute $\Gamma_i$ explicitly and from that get the desired rational expression of $x_i$ in 
  $x_1,\ldots,x_k$ over $K$, by dividing $\Gamma_i$ through $f^\prime(x_1)\ldots f^\prime(x_k)$. 

  Our formula for $\widehat{f_i}$ also gives us a rational expression of $x_i$ in $x_1,\ldots,x_k$ over $K$.  
  Since $f_i(T_i)$ is the minimal polynomial of $x_i$ over $K_{i-1}=L$, the degree of $f_i(T_i)$ is $1$. 
  Hence also the degree of $\widehat{f_i}$ in the variable $T_i$ is $1$. Thus the equation 
  $\widehat{f_i}(x_1,\ldots,x_i)=0$ is a a rational expression of $x_i$ in $x_1,\ldots,x_k$ over $K$. 
  We will exploit this idea in Section \ref{secgal} for special case of particular historical interest
  (concerning a theorem of Galois). In Section \ref{secexa}, we will present some examples of the same 
  special case. However, at that point we will not mention again that $\widehat{f_i}$
  serves as a polynomial expression of $x_i$ in $x_1,\ldots,x_k$ over $K$, 
  since readers will notice themselves. Both in Section \ref{secgal} and in Section \ref{secexa}
  we will work over the ground field $K=\Q$. 
\end{remark}

\begin{remark}
  Note also the structural similarity of the formula \eqref{mckay} with the sum 
  in \eqref{lag}. Nevertheless, the analogy does not go very far, since the $n$-tuple 
  $(\Gamma_1,\ldots,\Gamma_n)$ is a characteristic function of the set 
  $\{(\sigma(x_1),\ldots,\sigma(x_n));\sigma\in G\}$, whereas the polynomials $\widehat{f_i}$ all vanish 
  on precisely this set, as we shall discuss in the next section. 
\end{remark}

\begin{remark}
  Stauduhar and McKay (see \cite{mckay}) wrote that an effective method for the isolation of a generating family 
  of relations between the roots of the polynomial $f$ must be at least as powerful as a method for 
  finding the Galois group of $f$. Considering our interpolation formula, we see that is sufficient 
  to know the Galois group of $f$ in order to know a generating family of relations between the roots of $f$. 
  Nonetheless, we do need the Galois group in order to compute the generators $\widehat{f_i}$ 
  with our interpolation formula. Thus we answered the question,
  ``How can one do computations in the splitting field of a polynomial?''
  only up to an answer of the question, ``What is the Galois group of $f$?''

  It generally a difficult problem to compute the Galois group of a polynomial. Today's computer algebra 
  systems compute the Galois group using table-based methods for rational polynomials of degree $\leq23$. 
  For polynomials of arbitrarily high degree, no algorithm is known. In \cite{anai}, 
  the authors present a method of computing the Galois group of a polynomial which is based on 
  the factorisation of $f$ over $K_i$, $i=1,\ldots,n$. However, for an arbitrary ground field $K$, 
  there is no factorisation algorithm. Conversely, if we had a Gr\"obner basis of the relation ideal $I$, 
  from which source ever, it would be fairly easy to compute the Galois group of $f$. This is shown in \cite{anai}. 

  Thus we see that an effective method of isolating of a generating family 
  of relations between the roots of the polynomial $f$ is actually equivalent to a method of 
  finding the Galois group of $f$. In other words, the two questions,
  ``How can one do computations in the splitting field of a polynomial?'' and,
  ``What is the Galois group of $f$?'' are equivalent. 
\end{remark}

\section{Connection with the Buchberger--M\"oller algorithm}\label{secmoe}

At this point we explain how our work is far more deeply related to previous research. 
We show that our interpolation formula is a consequence of the Buchberger--M\"{o}ller algorithm. 
A very good introduction to the Buchberger--M\"{o}ller algorithm, 
besides the original article \cite{buchbergermoeller}, is \cite{rob}. 
For further reading, \cite{abbott} or \cite{alonso} are likewise very valuable. 
However, I take \cite{leo} as the best introduction. 

Let us briefly describe the situation for which the algorithm is designed. 
Let $F$ be a field and $\A^n(F)$ the affine $n$-space over $F$, 
and let $X$ be a finite set of points in $\A^n(F)$. As before, let $T=(T_1,\ldots,T_n)$ be indeterminates, 
this time over $F$. Let $I(X)$ be the ideal of $F[T]$ defining $X$, 
that is, $I(X)=\{g\in F[T];g(x)=0\text{ for all }x\in X\}$. The goal is to find a reduced 
Gr\"obner basis of the ideal $I(X)$ for a given monomial ordering of $F[T]$. 
The main objects that are used in the Buchberger--M\"{o}ller algorithm are the following:

Let $\deg(I(X))$ denote the subset of $\N^n$ that consists of all degrees of nonzero elements of $I(X)$, 
define $\Oh(X)=\N^n\setminus \deg(I(X))$. Since $I(X)$ is a zero-dimensional ideal, $\Oh(X)$ is finite. 
In particular, the residue classes of the monomials $T^\beta$, where $\beta$ runs through $\Oh(X)$, form a
basis of the $F$-vector space $F[T]/I(X)$. Furthermore, there exists a unique set of polynomials $h_x$, 
indexed by $x\in X$, such that $\deg(h_x)\in\Oh(X)$ for all $x\in X$, 
and such that $h_x(x)=1$ and $h_x(x^\prime)=0$ if $x^\prime\neq x$, for all $x,x^\prime\in X$. 
These polynomials are called the \textit{separators} of $X$; they can be used for solving a 
multivariate interpolation problem, namely: Given a collection of values $b_x\in F$, for every $x\in X$, 
there is a unique polynomial $b\in F[T]$ with support in $\Oh(X)$ such that $b(x)=b_x$ for all $x\in X$. Clearly
\begin{equation*}
  b=\sum_{x\in X}b_xh_x\,.
\end{equation*}
Finally, let $C(X)$ denote set of elements $T^\alpha$ of $\deg(I(X))$ such that $T^\beta$ lies in $\Oh(X)$ 
whenever $\alpha\neq\beta$ and $T^\beta$ divides $T^\alpha$. Think of $C(X)$ as the corners of $\deg(I(X))$. 
Since $\Oh(X)$ is finite, $C(X)$ is also finite, say $C(X)=\{T^{\alpha_1},\ldots,T^{\alpha_k}\}$. 

Were all these data given, it would be easy to compute the Gr\"obner basis: For $i=1,\ldots,k$, define 
\begin{equation}\label{basisels}
  g_i=T^{\alpha_i}-\sum_{x\in X}T^{\alpha_i}(x)h_x\,,
\end{equation}
then $\{g_1,\ldots,g_k\}$ is a Gr\"obner basis of $I(X)$. 
The problem is that, if given a finite set $X\subset\A^n(F)$, we do not a priori know the set $\Oh(X)$ 
nor the separators of $X$. It is the achievement of the Buchberger--M\"oller algorithm to compute these data, 
along with the Gr\"obner basis of $I(X)$. 

Now let us take a look at the results of the last two sections in the light of these ideas. 
Theorem \ref{procoe} says that in order to look for the Gr\"obner basis of the relation ideal $I$
(for the lexicographic ordering) of the polynomial $f$, we just need to look for polynomials 
$\widehat{f_i}$ vanishing at $(\sigma(x_1),\ldots,\sigma(x_i))$ for all $\sigma\in G$, 
and having the additional property that their degree looks as in \eqref{fhat}. What about this degree? 
It is easy to see that $\deg(I)$, the set of degrees of nonzero elements of $I$, 
equals $\cup_{i=1}^n(d_ie_i+\N^n)$. Here $d_i$ is the natural number as defined in Section \ref{secgen}, 
and $e_i$ is the $i$-th standard basis element of $\N^n$. Thus $\Oh(X)$, which we define to be 
$\N^n\setminus\deg(I(X))$, equals $\{\beta\in\N^n;\beta_i< d_i\text{ for all }i\}$. 
And the corners of $\deg(I)$ are the points $d_ie_i$, for $i=1,\ldots,n$. 
In particular, the degree of the polynomials in \eqref{lag} looks just like the degree of the polynomials in 
\eqref{basisels}. There is a still deeper similarity. Obviously, the polynomial $h_x$ 
(where $x$ is some point at which the ideal in question should vanish) in \eqref{basisels}, 
corresponds to the polynomial $\prod_{y_1\in B{(\rho,1)}}\frac{T_1-y_1}{\rho(x_1)-y_1}\ldots
\prod_{y_i\in B{(\rho,i)}}\frac{T_i-y_i}{\rho(x_i)-y_i}$ 
(which vanishes in $(\rho(x_1),\ldots,\rho(x_n))$) in \eqref{lag}. 

Let us formulate this observation more precisely. We prove the following theorem, 
in which we solve the problem of finding $I(X)$ as above, for $X$ of a particular shape. 

\begin{thm}\label{buch}
  Let $X\subset\A^n(F)$ be a finite set on which the group $G$ acts freely and transitively, 
  permuting the coordinates of the elements of $X$, i.e., $\sigma(x_i)=x_{\sigma(i)}$ 
  for all $\sigma\in G$ and for all $x=(x_1,\ldots,x_n)$ in $X$. Fix one $x\in X$. 
  Define $G_0=G$ and a chain of subgroups $G_i=\{\sigma\in G;\sigma(x)_j=x_j, j\leq i\}$ 
  for $i=1,\ldots,n$. Define $d_i=[G_{i-1}:G_i]$ and 
  $B{(\rho,i)}=\{\sigma(x)_i;\sigma\in G,\sigma(x)_j=\rho(x)_j, j\leq i-1\}\setminus\{\rho(x)_i\}$
  for $\rho\in G$. Then the polynomials 
  \begin{equation*}
    h_{\rho(x)}=\prod_{y_1\in B{(\rho,1)}}\frac{T_1-y_1}{\rho(x)_1-y_1}\ldots
    \prod_{y_n\in B{(\rho,n)}}\frac{T_n-y_n}{\rho(x)_n-y_n}\,,
  \end{equation*}
  for all $\rho\in G$, are the separators of $X$. Furthermore, the polynomials 
  \begin{equation*}
    g_i=T_i^{d_i}-\sum_{\rho\in G}\rho(x)_i^{d_i}h_{\rho(x)}\,,
  \end{equation*}
  for $i=1,\ldots,n$, form a Gr\"obner basis of $I(X)$ for the lexicographic ordering with 
  $T_1<\ldots<T_n$ on $F[T]$. In fact, we have 
  \begin{equation}\label{laggen}
    g_i=T_i^{d_i}-\sum_{\rho\in G/G_i}\rho(x)_i^{d_i}\prod_{y_1\in B{(\rho,1)}}\frac{T_1-y_1}{\rho(x)_1-y_1}\ldots
    \prod_{y_i\in B{(\rho,i)}}\frac{T_i-y_i}{\rho(x)_i-y_i}\,.
  \end{equation}
\end{thm}

\begin{pf}
  Take an arbitrary $\sigma\in G$. In the sum defining $g_i$, we distinguish the cases $\rho=\sigma$, 
  yielding $h_{\rho(x)}(\sigma(x))=1$, and $\rho\neq\sigma$, yielding $h_{\rho(x)}(\sigma(x))=0$. 
  Collecting terms, we get $g_i(\sigma(x))=0$ for all $\sigma\in G$. Thus $g_i$ lies in $I(X)$, 
  for all $i=1,\ldots,n$. 
  Define $I^\prime\subset I(X)$ to be the ideal of $F[T]$ generated by $g_i$, $i=1,\ldots,n$. This generating set is  
  triangular in the sense of \cite{aub}, thus a Gr\"obner basis of $I^\prime$. We conclude that
  the dimension of $F[T]/I^\prime$ as a vector space over $F$ equals $d_1\ldots d_n=|G|$. 

  On the other hand, for $\rho$ in $G$, let $I(\rho(x))$ be the vanishing ideal of the point $\rho(x)$ in $F[T]$. 
  Then $I(\rho(x))$ is generated by $T_i-\rho(x)_i$, $i=1,\ldots,n$. Clearly the quotient $F[T]/I(\rho(x))$ 
  is a one-dimensional vector space over $F$. For $\rho\neq\rho^\prime$, the ideals $I(\rho(x))$ and 
  $I(\rho^\prime(x))$ are coprime, and we clearly have $I(X)=\prod_{\rho\in G}I(\rho(x))$. 
  Thus the Chinese remainder theorem yields
  \begin{equation*}
    F[T]/I(X)=\oplus_{\rho\in G}F[T]/I(\rho(x))\,.
  \end{equation*}
  Since the quotient $F[T]/I(X)$ has the same $F$-dimension as the quotient $F[T]/I^\prime$ and  
  $I^\prime\subset I(X)$, it follows that $I^\prime=I(X)$. 

  In particular, $g_i$ are a Gr\"obner basis of $I(X)$, from which we see that
  $\Oh(X)=\{(k_1,\ldots,k_n);0\leq k_i\leq d_i-1\}$. Hence the polynomials $h_{\rho(x)}$ are supported in 
  $\Oh(X)$ and therefore are indeed the separators of $X$.

  To prove the last assertion, we split up the sum defining $g_i$ in the following way:
  \begin{equation}\label{better}
  \begin{split}
    &\sum_{\rho\in G}\rho(x)_i^{d_i}\prod_{y_1\in B{(\rho,1)}}\frac{T_1-y_1}{\rho(x)_1-y_1}\ldots
    \prod_{y_n\in B{(\rho,n)}}\frac{T_n-y_n}{\rho(x)_n-y_n}=\\
    &\sum_{\rho\in G/G_i}\rho(x)_i^{d_i}\prod_{y_1\in B{(\rho,1)}}\frac{T_1-y_1}{\rho(x)_1-y_1}\ldots
    \prod_{y_i\in B{(\rho,i)}}\frac{T_i-y_i}{\rho(x)_i-y_i}\\
    &\cdot\sum_{\tau\in G_i}\prod_{y_{i+1}\in B{(\rho\tau,i+1)}}\frac{T_{i+1}-y_{i+1}}{\rho\tau(x)_{i+1}-y_{i+1}}
    \ldots\prod_{y_n\in B{(\rho\tau,n)}}\frac{T_n-y_n}{\rho\tau(x)_n-y_n}\,.
  \end{split}
  \end{equation}
  $G_i$ operates freely and transitively on 
  $X(\rho,i) =\{(\rho\tau(x_{i+1}),\ldots,\rho\tau(x_n));\tau\in G_i\}$, thus by the first part of the theorem, 
  the factor in the last line of \eqref{better} is the unique polynomial with support in $\Oh(X(\rho,i))$ 
  and value $1$ at all points of $X(\rho,i)$. Clearly this polynomial is $1$. 
\end{pf}

\begin{remark}
  Note that the free and transitive operation of $G$ on $X$ is necessary in order to have equal $F$-dimensions
  of $F[T]/I(X)$ and $F[T]/I^\prime$, and that the operation as permutation of the coordinates is 
  necessary in order to define $d_i$ since this requirement guarantees that $G_i$ is a subgroup of $G_{i-1}$. 
\end{remark}

\begin{remark}
Of course this theorem applies to the situation where $F=L$ and 
$X=\{(\sigma(x_1),\ldots,\sigma(x_n));\sigma\in G\}$. 
The result for this special case is the interpolation formula \eqref{lag}. 
Now it seems awkward that we proved Theorem \ref{thmlag} first and afterwards proved 
another theorem from which the first one follows. But after all, I feed justified in so doing. 
The reason is that the proof of Theorem \ref{thmlag} used field theoretic arguments
(the degree of $\widehat{f_i}$, and the continuations of an automorphism of $K_{i-1}|K_0$ to 
$K_i|K_0$), whereas in the course of this section, we saw that field theory is not at all necessary to 
prove the interpolation formula. All one needs is the right kind of operation of the group $G$ on the 
set $X$. From the field theoretic point of view, this is certainly astounding. 
Nevertheless, the fact that all polynomials $\widehat{f_i}$ are defined over $K$ and not over $F=L$
does need field theoretic reasoning, as we carried out in Section \ref{secgen}. 
\end{remark}

\begin{remark}
One can also prove an interpolation formula for $\widehat{f_i}$ 
(hence also for the more general situation that we discussed in this section) 
which is more in the spirit of Newton than Lagrange interpolation. 
This has been done in \cite{led}. However, the generalisation of the Newton formula to our multivariate situation
has the weakness of being rather messy and not leading to any deeper insight. This is why I do not present it here. 
\end{remark}

\section{On a theorem of Galois}\label{secgal}

The following theorem is equivalent to the theorem mentioned in the introduction and due to Evariste Galois 
(see e.g. \cite{hup}, Satz 21.4, p.612, or \cite{cohn}, Theorem 11.7, p.119):

\begin{thm}[Galois]\label{thmgal}
  Let $f\in\Q[Z]$ be an irreducible polynomial of degree $p$, where $p$ is a prime number. 
  Let $L$ be the splitting field of $f$ and $x_1,\ldots,x_p$ the zeros of $f$ in $L$. 
  Then the following statements are equivalent:
  \begin{enumerate}
  \item[(i)]$f$ is solvable by radicals.
  \item[(ii)]$L=\Q(x_i,x_j)$ for all $i,j\in\{1,\ldots,p\}$ such that $i\neq j$.
  \end{enumerate}
\end{thm}

\begin{remark}
  In the terminology of \cite{mckay}, the theorem reads as follows: A rational irreducible polynomial $f$
  is solvable if and only if any two zeros of $f$ form a minimal strong base of $f$ over $\Q$. 
\end{remark}

We would like to apply the results of Sections \ref{secgen} and \ref{secint} in order to give an explicit formula 
for the polynomial dependence of $x_k$ from $x_i$ and $x_j$. By the theorem, any other zero $x_k$ lies in 
$\Q(x_i,x_j)$. We choose a numbering of the zeros such that $x_i=x_1$, $x_j=x_2$, $x_k=x_3$. We determine the 
minimal polynomial $\widehat{f_3}$ of $x_3$ over $K_2$. It has degree 1, thus can be written 
$\widehat{f_3}=T_3-P(x_1,x_2)$ for a suitable $P\in K[T_1,T_2]$. We evaluate $\widehat{f_3}$ at $x_3$ and obtain 
an equation $x_3=P(x_1,x_2)$. This is the rational polynomial in two zeros, whose existence is claimed in the 
theorem. Recalling the precise form of $\widehat{f_3}$ in \ref{lag}, we obtain
\begin{equation*}
  x_3=\sum_{\rho \in G/G_3}\rho(x_3)\prod_{y_1\in B{(\rho,1)}}\frac{x_1-y_1}{\rho(x_1)-y_1}
  \prod_{y_2\in B{(\rho,2)}}\frac{x_2-y_2}{\rho(x_2)-y_2}\,.
\end{equation*}
Note that it is not clear a priori that the right hand side of this formula is a rational polynomial in $x_1$ and $x_2$!

Let us now assume that $G$ is 2-fold sharply transitive, that is, that $G$ operates freely and transitively  
on the set of pairs $\{(x_i,x_j);i\neq j\}$. In this case, the formula takes a particularly neat shape:
$G_3$ is trivial and $G/G_3$ can be identified with the set of all pairs $(x,y)$ of distinct zeros. 
Each such pair can be written as $(x,y)=(\rho(x_1),\rho(x_2))$ for a unique $\rho\in G$. We define 
$z(x,y)=\rho(x_3)$ and note that 
$B{(\rho,1)}=\{x_1,\ldots,x_p\}\setminus\{x\}$ and $B{(\rho,2)}=\{x_1,\ldots,x_p\}\setminus\{x,y\}$. Thus we obtain
\begin{equation*}
  x_3=\sum_{x\neq y\in\{x_1,\ldots,x_p\}}\prod_{\substack{\xi\in\{x_1,\ldots,x_p\}\\\xi\neq x}}\frac{x_1-\xi}{x-\xi}
  \prod_{\substack{\eta\in\{x_1,\ldots,x_p\}\\\eta\neq x,y}}z(x,y)\frac{x_2-\eta}{y-\eta}\,.
\end{equation*}

\section{Numerical computation of the generators}\label{secnum}

In this section we always assume $K=\Q$, and work out the technicalities that will enable us to 
practically compute the generators of the relation ideal for some given rational polynomials in Section \ref{secexa}. 
So let $f$ over $\Q$ be given. We will work with complex and $p$-adic approximations of the zeros 
of $f$ in order to construct the generators of the relation ideal. This task requires the knowledge of an integer 
$\Delta_i$ such that $\Delta_i\widehat{f_i}$ has integer coefficients (Theorem \ref{prointeger}). 
Further, we need an upper bound for the absolute values of these coefficients (Theorem \ref{probound}). The final 
result is formulated in Theorem \ref{procom}. 

Let $\gamma$ be a rational integer such that all the products $\gamma x_j$, $j=1,\ldots,n$, are algebraic integers. 
We denote the discriminant of $f$ by 
\begin{equation*}
  d(f)=\prod_{1\leq r<s\leq n}(x_r-x_s)^2\,.
\end{equation*}
The ceiling function is always denoted by $\lceil\,\,\,\rceil$ and the floor function by $\lfloor\,\,\,\rfloor$. 

\begin{thm}\label{prointeger}
  For $i=1,\ldots,n$, the rational integer 
  \begin{equation*}
    \Delta_i=\gamma^{n(n-1)\lceil\frac{i}{2}\rceil+d_i}d(f)^{\lceil\frac{i}{2}\rceil}
  \end{equation*}
  has the property that $\Delta_i\widehat{f_i}$ lies in $\Z[T_1,\ldots,T_i]$. 
\end{thm}

\begin{pf}
  First observe that $\gamma^{n(n-1)}d(f)$ lies in $\Z$,, since its factors $\gamma(x_r-x_s)$ lie in $\Z$. 
  Thus also all $\Delta_i$ lie in $\Z$. In Section \ref{secint} we proved the interpolation formula \eqref{lag}. 
  We multiply this equation by $d(f)^{\lceil\frac{i}{2}\rceil}$. On the right hand side, 
  we have the summand $d(f)^{\lceil\frac{i}{2}\rceil}T_i^{d_i}$ and for all $\rho\in G/G_i$ a summand of the form
  \begin{equation*}
    \prod_{1\leq r<s\leq n}(x_r-x_s)^{2\lceil\frac{i}{2}\rceil}\rho(x_i)^{d_i}
    \prod_{y_1\in B{(\rho,1)}}\frac{T_1-y_1}{\rho(x_1)-y_1}\ldots\prod_{y_i\in B{(\rho,i)}}
    \frac{T_i-y_i}{\rho(x_i)-y_i}\,.
  \end{equation*}
  Clearly the denominators in this fraction are cancelled if multiplied by suitable factors of the leftmost product. 
  After reduction there remain $n(n-1)\lceil\frac{i}{2}\rceil+d_i$ factors in the product. 
  All of them have the property that 
  if multiplied by $\gamma$ they are polynomials in $L[T_1,\ldots,T_i]$ 
  with algebraic integers as coefficients. Thus also the product $\Delta_i\widehat{f_i}$ is a polynomial in 
  $L[T_1,\ldots,T_i]$ with algebraic integers as coefficients. On the other hand, this product lies in 
  $\Q[T_1,\ldots,T_i]$, hence it lies in $\Z[T_1,\ldots,T_i]$, as claimed.
\end{pf}

For the time being, let $x_1,\ldots,x_n\in\C$ denote the complex zeros of $f$ and $\abs{\,\,\,\,}$ 
denote the usual absolute value in $\C$. 

\begin{thm}\label{probound}
  Let $D,M\in \R_{>0}$ be such that $M>\max\{\abs{x_r}\}$ and $D>\max\{\abs{x_r-x_s};x_r\neq x_s\}$. 
  Then the absolute value of $\Delta_ib_{i,k_1,\ldots,k_i}$ is bounded by 
  \begin{equation}\label{prodi}
    \gamma^{n(n-1)\lceil\frac{i}{2}\rceil+d_i}\binom{d_1-1}{k_1-1}\ldots\binom{d_i-1}{k_i-1}
    M^{k_1+\ldots+k_i-i+d_i}D^{n(n-1)\lceil\frac{i}{2}\rceil-d_1-\ldots-d_i+i}\,.
  \end{equation}
\end{thm}

\begin{pf}
  We evaluate the formula \eqref{lag} for $\widehat{f_i}$ at the complex zeros and multiply the result by $\Delta_i$. 
  As in the proof of Theorem \ref{prointeger} we cancel the denominators occurring in the sum by factors of the 
  product $d(f)^{\lceil\frac{i}{2}\rceil}$ which occur in $\Delta_i$. In the remaining product we have 
  $n(n-1)\lceil\frac{i}{2}\rceil-d_1-\ldots-d_i+i$ factors of the type $(\xi_r-\xi_s)$ left. 
  The absolute value of these is bounded by $D$. Further, $M$ is an upper bound for $\rho(\xi_i)$. 
  Finally, it is easy to check that for $j=1,\ldots,i$, the absolute value of the coefficient of the polynomial 
  $\prod_{y_j\in B{(\rho,j)}}(T_j-y_j)$ at $T_j^{d_j-k_j}$ is bounded by $\binom{d_j-1}{k_j-1}M^{k_j-1}$. 
  We obtain the result by collecting factors.
\end{pf}

We fix an integer $c$ such that $cf$ lies in $\Z[Z]$. Let $p$ be a prime number such that the 
polynomial $\overline{cf}\in (\Z/p\Z)[Z]$ (the reduction of $cf$ modulo $p$) splits into $n=\deg(f)$ 
disjoint linear factors over $\Z/p\Z$. By Hensel's lemma, we can lift these zeros to zeros of $cf$ in $\Q_p$. 
The polynomial $cf$ also 
splits into $n$ disjoint linear factors over $\Z/p^e\Z$, for all integers $e\leq 1$. In this process, if $e<k$, 
the zeros in $\Z/p^e\Z$ are obtained from the zeros in $\Z/p^k\Z$ by reduction modulo $p^e$. 
We call the zeros in $\Z/p^e\Z$ the $e$-th $p$-adic approximations of the zeros lying in $\Q_p$. 

The existence of a prime with the property mentioned above follows from Chebotarev's density theorem, 
see e.g. \cite{neu}. The density of these primes equals $1/|G|$; thus if the Galois group is large, it might become
difficult to find such a prime. However, in this case case one could also choose a prime $p$ without the 
property mentioned above and work with a suitable algebraic extension of $\Q_p$ instead of with $\Q_p$ itself. 
We will not discuss this topic any further since for practical purposes, the Galois group 
never gets too big, and it is no problem at all to find a suitable $p$. 

For the forthcoming discussion, we let $G$ operate on the $p$-adic approximations of the zeros in the obvious way.  
We will need $p$-adic approximations of $d(f)$, $B{(\rho,i)}$, $\widehat{f_i}$ and $\Delta_i$. The approximations
are defined by the same formulas as the original objects, but with each zero replaced by the respective approximation. 
Now we can specify exponents $e_i$ such that from the knowledge of $e_i$-th $p$-adic approximations of the zeros 
of $f$ we can compute $\widehat{f_i}$. 

\begin{thm}\label{procom}
  For $i=1,\ldots,n$ the following holds: Let $\lambda_i$ be the maximum of $\abs{\Delta_i}$ and the products 
  \eqref{prodi}, for all $k_j=1,\ldots,d_j$. Define $e_i=\lfloor\frac{\log(2\lambda_i-1)}{\log(p)}\rfloor+1$. 
  We view the $e_i$-th $p$-adic approximation of $\Delta_i\widehat{f_i}$ as a polynomial in $\Z[T_1,\ldots,T_i]$ 
  by using the system of representatives $\{-\frac{p^{e_i}-1}{2},\ldots,\frac{p^{e_i}-1}{2}\}$ of $\Z/p^{e_i}\Z$. 
  Then this polynomial coincides with $\Delta_i\widehat{f_i}$. 
\end{thm}

\begin{pf}
  Let $\beta_{i,k_1,\ldots,k_i}\in\{-\frac{p^{e_i}-1}{2},\ldots,\frac{p^{e_i}-1}{2}\}$ be the coefficients of the 
  approximate $\Delta_i\widehat{f_i}$. 
  Then we can find $\mu_{i,k_1,\ldots,k_i}\in\Z$ such that 
  $b_{i,k_1,\ldots,k_i}=\beta_{i,k_1,\ldots,k_i}+\mu_{i,k_1,\ldots,k_i}p^{e_i}$.
  Now if $\mu_{i,k_1,\ldots,k_i}$ were not zero, we would have $\abs{b_{i,k_1,\ldots,k_i}}\geq(p^{e_i}+1)/2$.
  On the other hand, by definition of $e_i$ we have $e_i>\log(2\lambda_i-1)/\log(p)$, from which we deduce 
  $\lambda_i<(p^{e_i}+1)/2$. We have assumed $\abs{b_{i,k_1,\ldots,k_i}}<\lambda_i$, hence 
  $\abs{b_{i,k_1,\ldots,k_i}}<(p^{e_i}+1)/2$, a contradiction. Thus $\mu_{i,k_1,\ldots,k_i}=0$. 
\end{pf}

\section{Examples}\label{secexa}

In this section we present some examples treated with the methods developed in Section \ref{secnum}. 
We used KANT (\cite{KANT}) and GAP (\cite{GAP4}) for all computations. 
GAP was used for the computation of $\widehat{f_i}$ by the 
interpolation formula, since this formula requires computations in permutation groups, 
the purpose for which GAP is designed. KANT was used for computing $p$-adic approximations of the zeros and 
for checking the Galois group of the sample polynomials. KANT can compute the Galois group of 
irreducible polynomials of degree $\leq23$ over $\Q$. We point out that KANT is not the only computer algebra system
capable of computing the Galois group. PARI, MAPLE or MAGMA, just to mention a few, can also do this, 
The reason why we decided to use KANT and GAP is that these systems are freely available. 
Until recently, KANT and the other systems computed the Galois group of a polynomial
only up to conjugacy (by a table-based method and not a direct method, in the terminology of \cite{anai}). 
This used to cause quite a few problems in donating examples, 
because the interpolation formula \eqref{lag} requires knowledge of $G$ as acting on the zeros of $f$. 
Meanwhile the situation has greatly improved, since KANT features some commands that allow us 
to compute the Galois group of $f$ acting on approximations of the zeros of $f$. 
Obviously KANT is now using a direct method, in the terminology of \cite{anai}, in doing this computation. 

Here are the examples. For various irreducible separable polynomials over $\Q$, we give the following data: 
The Galois group by all the names it bears in the database created by J\"urgen Kl\"uners and Gunter Malle, to be found at
{\verb http://www.mathematik.uni-kassel.de/~klueners/minimum/minimum.html } 
(we have embedded it into the appropriate symmetric group $S_n$ by using the GAP command {\verb TransitiveGroup }), 
the discriminant of $f$, a prime $p$ as in Section \ref{secnum}, the approximation
in $\Z/p\Z$ of the zeros $x=(x_1,\ldots,x_n)$ to which corresponds the chosen embedding of the Galois group into $S_n$, 
and the generators $\widehat{f_i}\,,\,i=1,\ldots,n$ of the relation ideal. Note that $\widehat{f_1}=f(T_1)$, 
so in each example we compute the relation ideal of $\widehat{f_1}$. The monomials of $\widehat{f_i}$ 
are ordered lexicographically with $T_1<\ldots<T_n$. The polynomials $\widehat{f_1}$ 
were taken from the Kl\"uners--Malle database mentioned above. 

\begin{example}\label{gal1}
$G = 5T1 = C(5) = 5$, $p=23$, $x=(19, 9, 13, 17, 12)$, $d(f)=14641$, \\
$\widehat{f_1} = T_{1}^{5}-T_{1}^{4}-4T_{1}^{3}+3T_{1}^{2}+3T_{1}-1$, \\
$\widehat{f_2} = T_{2}+T_{1}^{3}-3T_{1}+2$, \\
$\widehat{f_3} = T_{3}-T_{1}^{4}+2T_{1}^{3}-T_{1}^{2}-T_{1}+1$, \\
$\widehat{f_4} = T_{4}+T_{1}^{4}-T_{1}^{3}-2T_{1}^{2}+2T_{1}-1$, \\
$\widehat{f_5} = T_{5}-2T_{1}^{3}+3T_{1}^{2}+2T_{1}-2$.
\end{example}

\begin{example}\label{gal2}
$G = 5T3 = F(5) = 5 : 4$, $p=191$, $x=(70, 61, 33, 11, 17)$, $d(f)=35152$, \\
$\widehat{f_1} = T_{1}^{5}-T_{1}^{4}+2T_{1}^{3}-4T_{1}^{2}+T_{1}-1$, \\
$\widehat{f_2} = T_{2}^{4}+2/13T_{1}^{4}T_{2}^{3}+8/13T_{1}^{3}T_{2}^{3}-30/13T_{1}^{2}T_{2}^{3}
+4/13T_{1}T_{2}^{3}-14/13T_{2}^{3}-7/26T_{1}^{4}T_{2}^{2}-15/13T_{1}^{3}T_{2}^{2}
+12/13T_{1}^{2}T_{2}^{2}+2/13T_{1}T_{2}^{2}+1/2T_{2}^{2}+11/13T_{1}^{4}T_{2}-2/13T_{1}^{3}T_{2}
-32/13T_{1}T_{2}-8/13T_{2}-3/26T_{1}^{4}-6/13T_{1}^{3}+3/13T_{1}^{2}-3/13T_{1}+21/26$, \\
$\widehat{f_3} = T_{3}-8/13T_{1}^{4}T_{2}^{3}+2/13T_{1}^{3}T_{2}^{3}-11/13T_{1}^{2}T_{2}^{3}
+3/13T_{1}T_{2}^{3}-8/13T_{2}^{3}+15/13T_{1}^{4}T_{2}^{2}-8/13T_{1}^{3}T_{2}^{2}
+11/13T_{1}^{2}T_{2}^{2}-21/13T_{1}T_{2}^{2}+10/13T_{2}^{2}-5/13T_{1}^{4}T_{2}
+10/13T_{1}^{3}T_{2}-36/13T_{1}^{2}T_{2}+32/13T_{1}T_{2}-23/13T_{2}+9/13T_{1}^{4}
-9/13T_{1}^{3}+12/13T_{1}^{2}-19/13T_{1}-2/13$, \\
$\widehat{f_4} = T_{4}+11/26T_{1}^{4}T_{2}^{3}-9/13T_{1}^{3}T_{2}^{3}+14/13T_{1}^{2}T_{2}^{3}
-9/13T_{1}T_{2}^{3}+3/26T_{2}^{3}-1/2T_{1}^{4}T_{2}^{2}-4/13T_{1}^{3}T_{2}^{2}
+20/13T_{1}^{2}T_{2}^{2}-3/13T_{1}T_{2}^{2}+25/26T_{2}^{2}-31/26T_{1}^{4}T_{2}
+45/13T_{1}^{3}T_{2}-62/13T_{1}^{2}T_{2}+46/13T_{1}T_{2}-5/2T_{2}+33/26T_{1}^{4}
-15/13T_{1}^{3}+21/13T_{1}^{2}-31/13T_{1}-1/26$, \\
$\widehat{f_5} = T_{5}+3/26T_{1}^{4}T_{2}^{3}-3/13T_{1}^{3}T_{2}^{3}-23/13T_{1}^{2}T_{2}^{3}
+6/13T_{1}T_{2}^{3}-7/26T_{2}^{3}+19/26T_{1}^{4}T_{2}^{2}-16/13T_{1}^{3}T_{2}^{2}
+17/13T_{1}^{2}T_{2}^{2}-2/13T_{1}T_{2}^{2}+3/26T_{2}^{2}-19/26T_{1}^{4}T_{2}
+9/13T_{1}^{3}T_{2}-21/13T_{1}^{2}T_{2}-21/26T_{2}+17/26T_{1}^{4}-3/13T_{1}^{2}-2/13T_{1}+9/26$.
\end{example}

\begin{example}
$G = 6T3 = D(6) = S(3)[x]2$, 
$p=211$, $x=(203, 22, 127, 99, 92, 90)$, $d(f)=-14283$, \\
$\widehat{f_1} = T_{1}^{6}+T_{1}^{4}-2T_{1}^{3}+T_{1}^{2}-T_{1}+1$, \\
$\widehat{f_2} = T_{2}^{2}+9/23T_{1}^{5}T_{2}+27/23T_{1}^{4}T_{2}-15/23T_{1}^{3}T_{2}+22/23T_{1}^{2}T_{2}
-30/23T_{1}T_{2}+6/23T_{2}-27/23T_{1}^{5}+6/23T_{1}^{4}-9/23T_{1}^{3}+21/23T_{1}^{2}+3/23T_{1}-16/23$, \\
$\widehat{f_3} = T_{3}-4/3T_{1}^{4}T_{2}+T_{1}^{3}T_{2}-2/3T_{1}^{2}T_{2}+1/3T_{1}T_{2}-T_{2}$, \\
$\widehat{f_4} = T_{4}+12/23T_{1}^{5}T_{2}+9/23T_{1}^{4}T_{2}+33/23T_{1}^{3}T_{2}-21/23T_{1}^{2}T_{2}
+45/23T_{1}T_{2}-27/23T_{2}-9/23T_{1}^{5}-50/23T_{1}^{4}-12/23T_{1}^{3}-16/23T_{1}^{2}+44/23T_{1}-6/23$, \\
$\widehat{f_5} = T_{5}-3T_{1}^{5}T_{2}-T_{1}^{3}T_{2}+3T_{1}^{2}T_{2}-T_{1}T_{2}+T_{2}-3T_{1}^{3}
+2T_{1}^{2}-3T_{1}+3$, \\
$\widehat{f_6} = T_{6}-2T_{1}^{4}T_{2}-T_{1}^{2}T_{2}+2T_{1}T_{2}-T_{1}^{5}-T_{1}^{3}+T_{1}^{2}-1$.
\end{example}

\begin{example}
$G = 6T5 = F_18(6) = [3^2]2 = 3wr2$, 
$p=151$, \\ $x=(114, 103, 46, 93, 25, 75)$, $d(f)=-9747$, \\
$\widehat{f_1} = T_{1}^{6}-3T_{1}^{5}+4T_{1}^{4}-2T_{1}^{3}+T_{1}^{2}-T_{1}+1$, \\
$\widehat{f_2} = T_{2}^{3}-31/19T_{1}^{5}T_{2}^{2}+4T_{1}^{4}T_{2}^{2}-113/19T_{1}^{3}T_{2}^{2}
+72/19T_{1}^{2}T_{2}^{2}-18/19T_{1}T_{2}^{2}-16/19T_{2}^{2}+22/19T_{1}^{5}T_{2}
-3T_{1}^{4}T_{2}+80/19T_{1}^{3}T_{2}-40/19T_{1}^{2}T_{2}+8/19T_{1}T_{2}
+T_{2}-3/19T_{1}^{5}+23/19T_{1}^{3}-64/19T_{1}^{2}+50/19T_{1}-12/19$, \\
$\widehat{f_3} = T_{3}-T_{1}^{5}T_{2}^{2}+5T_{1}^{3}T_{2}^{2}-6T_{1}^{2}T_{2}^{2}+3T_{1}T_{2}^{2}
+2/3T_{1}^{5}T_{2}+8/3T_{1}^{4}T_{2}-8T_{1}^{3}T_{2}+20/3T_{1}^{2}T_{2}-4/3T_{1}T_{2}
-2/3T_{2}+3T_{1}^{5}-14/3T_{1}^{4}+7T_{1}^{3}-16/3T_{1}^{2}+8/3T_{1}$, \\
$\widehat{f_4} = T_{4}+4/19T_{1}^{5}T_{2}^{2}-T_{1}^{4}T_{2}^{2}+3T_{1}^{3}T_{2}^{2}
-98/19T_{1}^{2}T_{2}^{2}+67/19T_{1}T_{2}^{2}-11/19T_{2}^{2}-41/19T_{1}^{5}T_{2}+6T_{1}^{4}T_{2}
-8T_{1}^{3}T_{2}+140/19T_{1}^{2}T_{2}-93/19T_{1}T_{2}+32/19T_{2}+35/19T_{1}^{5}-5T_{1}^{4}+7T_{1}^{3}
-107/19T_{1}^{2}+59/19T_{1}-6/19$, \\
$\widehat{f_5} = T_{5}-T_{1}^{5}T_{2}^{2}+5T_{1}^{4}T_{2}^{2}-10T_{1}^{3}T_{2}^{2}+10T_{1}^{2}T_{2}^{2}
-6T_{1}T_{2}^{2}+3T_{2}^{2}+4/3T_{1}^{5}T_{2}-22/3T_{1}^{4}T_{2}+14T_{1}^{3}T_{2}-12T_{1}^{2}T_{2}
+4T_{1}T_{2}-4/3T_{2}+4T_{1}^{5}-37/3T_{1}^{4}+55/3T_{1}^{3}-32/3T_{1}^{2}+T_{1}+11/3$, \\
$\widehat{f_6} = T_{6}+29/19T_{1}^{5}T_{2}^{2}-4T_{1}^{4}T_{2}^{2}+6T_{1}^{3}T_{2}^{2}
-74/19T_{1}^{2}T_{2}^{2}+25/19T_{1}T_{2}^{2}+20/19T_{2}^{2}-17/19T_{1}^{5}T_{2}+2T_{1}^{4}T_{2}
-4T_{1}^{3}T_{2}+8/19T_{1}^{2}T_{2}+43/19T_{1}T_{2}-72/19T_{2}-17/19T_{1}^{5}+2T_{1}^{4}-T_{1}^{3}
-11/19T_{1}^{2}+5/19T_{1}+23/19$.
\end{example}

\begin{example}\label{gal3}
$G = 7T1 = C(7) = 7$, 
$p=41$, $x=(6, 39, 35, 26, 21, 30, 8)$, $d(f)=171903939769$, \\
$\widehat{f_1} = T_{1}^{7}-T_{1}^{6}-12T_{1}^{5}+7T_{1}^{4}+28T_{1}^{3}-14T_{1}^{2}-9T_{1}-1$, \\
$\widehat{f_2} = T_{2}+3/17T_{1}^{6}-5/17T_{1}^{5}+11/17T_{1}^{4}-15/17T_{1}^{3}-4/17T_{1}-1/17$, \\
$\widehat{f_3} = T_{3}+4/17T_{1}^{6}-5/17T_{1}^{5}-31/17T_{1}^{4}+24/17T_{1}^{3}+5/17T_{1}^{2}+36/17T_{1}+9/17$, \\
$\widehat{f_4} = T_{4}-12/17T_{1}^{6}+14/17T_{1}^{5}+4T_{1}^{4}-33/17T_{1}^{3}-84/17T_{1}^{2}+12/17T_{1}+4/17$, \\
$\widehat{f_5} = T_{5}+5/17T_{1}^{6}-11/17T_{1}^{5}-40/17T_{1}^{4}-19/17T_{1}^{3}+88/17T_{1}^{2}-9/17T_{1}-7/17$, \\
$\widehat{f_6} = T_{6}+6/17T_{1}^{6}+14/17T_{1}^{5}-12/17T_{1}^{4}-10/17T_{1}^{3}+21/17T_{1}^{2}-24/17T_{1}-4/17$, \\
$\widehat{f_7} = T_{7}-6/17T_{1}^{6}-7/17T_{1}^{5}+4/17T_{1}^{4}+53/17T_{1}^{3}-30/17T_{1}^{2}-11/17T_{1}-1/17$.
\end{example}

\begin{example}\label{gal4}
$G = 7T3 = F_21(7) = 7:3$, 
$p=167$, $x=(71, 153, 6, 58, 23, 53, 137)$, $d(f)=1817487424$, \\
$\widehat{f_1} = T_{1}^{7}-8T_{1}^{5}-2T_{1}^{4}+16T_{1}^{3}+6T_{1}^{2}-6T_{1}-2$, \\
$\widehat{f_2} = T_{2}^{3}-53/73T_{1}^{6}T_{2}^{2}-89/73T_{1}^{5}T_{2}^{2}+369/73T_{1}^{4}T_{2}^{2}
+692/73T_{1}^{3}T_{2}^{2}-486/73T_{1}^{2}T_{2}^{2}-1284/73T_{1}T_{2}^{2}-446/73T_{2}^{2}
-142/73T_{1}^{6}T_{2}+1245/73T_{1}^{4}T_{2}+612/73T_{1}^{3}T_{2}-2054/73T_{1}^{2}T_{2}
-1550/73T_{1}T_{2}-228/73T_{2}-165/73T_{1}^{6}-12/73T_{1}^{5}+1137/73T_{1}^{4}+208/73T_{1}^{3}
-1922/73T_{1}^{2}-490/73T_{1}+472/73$, \\
$\widehat{f_3} = T_{3}+21/73T_{1}^{6}T_{2}^{2}+17/73T_{1}^{5}T_{2}^{2}-254/73T_{1}^{4}T_{2}^{2}
-184/73T_{1}^{3}T_{2}^{2}+429/73T_{1}^{2}T_{2}^{2}+222/73T_{1}T_{2}^{2}-82/73T_{2}^{2}
+103/73T_{1}^{6}T_{2}-29/73T_{1}^{5}T_{2}-656/73T_{1}^{4}T_{2}+79/73T_{1}^{3}T_{2}
+1051/73T_{1}^{2}T_{2}-170/73T_{1}T_{2}-516/73T_{2}+48/73T_{1}^{6}+9/73T_{1}^{5}
-235/73T_{1}^{4}+62/73T_{1}^{3}+472/73T_{1}^{2}-182/73T_{1}-398/73$, \\
$\widehat{f_4} = T_{4}+88/73T_{1}^{6}T_{2}^{2}+92/73T_{1}^{5}T_{2}^{2}-512/73T_{1}^{4}T_{2}^{2}
-850/73T_{1}^{3}T_{2}^{2}-1/73T_{1}^{2}T_{2}^{2}+497/73T_{1}T_{2}^{2}+151/73T_{2}^{2}
+104/73T_{1}^{6}T_{2}+235/73T_{1}^{5}T_{2}-671/73T_{1}^{4}T_{2}-1344/73T_{1}^{3}T_{2}
+195/73T_{1}^{2}T_{2}+824/73T_{1}T_{2}+208/73T_{2}+131/73T_{1}^{6}+79/73T_{1}^{5}
-557/73T_{1}^{4}-419/73T_{1}^{3}+404/73T_{1}^{2}+322/73T_{1}+6/73$, \\
$\widehat{f_5} = T_{5}-58/73T_{1}^{6}T_{2}^{2}-85/73T_{1}^{5}T_{2}^{2}+295/73T_{1}^{4}T_{2}^{2}
+446/73T_{1}^{3}T_{2}^{2}-328/73T_{1}^{2}T_{2}^{2}-525/73T_{1}T_{2}^{2}-152/73T_{2}^{2}
+7/73T_{1}^{6}T_{2}-51/73T_{1}^{5}T_{2}-44/73T_{1}^{4}T_{2}+531/73T_{1}^{3}T_{2}
+71/73T_{1}^{2}T_{2}-1146/73T_{1}T_{2}-532/73T_{2}+8/73T_{1}^{6}-19/73T_{1}^{5}
+62/73T_{1}^{4}+164/73T_{1}^{3}-118/73T_{1}^{2}-340/73T_{1}-298/73$, \\
$\widehat{f_6} = T_{6}+1/73T_{1}^{6}T_{2}^{2}+28/73T_{1}^{5}T_{2}^{2}+83/73T_{1}^{4}T_{2}^{2}
+43/73T_{1}^{3}T_{2}^{2}-31/73T_{1}^{2}T_{2}^{2}+21/73T_{1}T_{2}^{2}+12/73T_{2}^{2}
-75/73T_{1}^{6}T_{2}-187/73T_{1}^{5}T_{2}+291/73T_{1}^{4}T_{2}+1088/73T_{1}^{3}T_{2}
+384/73T_{1}^{2}T_{2}-449/73T_{1}T_{2}-170/73T_{2}-48/73T_{1}^{6}-159/73T_{1}^{5}
+218/73T_{1}^{4}+892/73T_{1}^{3}+281/73T_{1}^{2}-384/73T_{1}-138/73$, \\
$\widehat{f_7} = T_{7}-67/73T_{1}^{6}T_{2}^{2}-105/73T_{1}^{5}T_{2}^{2}+178/73T_{1}^{4}T_{2}^{2}
+418/73T_{1}^{3}T_{2}^{2}+239/73T_{1}^{2}T_{2}^{2}-80/73T_{1}T_{2}^{2}-82/73T_{2}^{2}
-19/73T_{1}^{6}T_{2}-42/73T_{1}^{5}T_{2}-46/73T_{1}^{4}T_{2}+2T_{1}^{3}T_{2}
+405/73T_{1}^{2}T_{2}+147/73T_{1}T_{2}-110/73T_{2}+93/73T_{1}^{6}+41/73T_{1}^{5}
-534/73T_{1}^{4}-234/73T_{1}^{3}+551/73T_{1}^{2}+114/73T_{1}-112/73$.
\end{example}

\begin{example}
$G = 8T20 = [2^3]4$, 
$p=337$, $x=(286, 282, 254, 175, 51, 55, 83, 162)$, $d(f)=4398046511104$, \\
$\widehat{f_1} = T_{1}^{8}-4T_{1}^{6}-6T_{1}^{4}+4T_{1}^{2}+1$, \\
$\widehat{f_2} = T_{2}^{2}+1/2T_{1}^{6}T_{2}-T_{1}^{4}T_{2}+1/2T_{1}^{2}T_{2}+1/2T_{1}^{5}-1/2T_{1}$, \\
$\widehat{f_3} = T_{3}^{2}-1/4T_{1}^{7}T_{2}T_{3}+T_{1}^{5}T_{2}T_{3}-1/4T_{1}^{3}T_{2}T_{3}
-2T_{1}^{6}T_{3}-7/4T_{1}^{4}T_{3}+T_{1}^{2}T_{3}+1/4T_{3}+7/4T_{1}^{4}T_{2}-T_{1}^{2}T_{2}
-1/4T_{2}-1/4T_{1}^{5}-T_{1}^{3}-1/4T_{1}$, \\
$\widehat{f_4} = T_{4}+1/2T_{1}^{6}T_{2}T_{3}+T_{1}^{4}T_{2}T_{3}+1/2T_{1}^{2}T_{2}T_{3}
-3/4T_{1}^{7}T_{3}+1/4T_{1}^{5}T_{3}+9/4T_{1}^{3}T_{3}+1/4T_{1}T_{3}-1/2T_{1}^{7}T_{2}
+1/2T_{1}^{5}T_{2}+3/2T_{1}^{3}T_{2}-1/2T_{1}T_{2}-1/4T_{1}^{6}+7/4T_{1}^{4}+15/4T_{1}^{2}+3/4$, \\
$\widehat{f_5} = T_{5}+1/4T_{1}^{6}T_{2}T_{3}-7/4T_{1}^{4}T_{2}T_{3}+5/4T_{1}^{2}T_{2}T_{3}
+1/4T_{2}T_{3}-1/4T_{1}^{7}T_{3}+3/4T_{1}^{5}T_{3}+11/4T_{1}^{3}T_{3}+3/4T_{1}T_{3}
+1/4T_{1}^{7}T_{2}-3/4T_{1}^{5}T_{2}-11/4T_{1}^{3}T_{2}+5/4T_{1}T_{2}+3/4T_{1}^{6}
-11/4T_{1}^{4}-21/4T_{1}^{2}+1/4$, \\
$\widehat{f_6} = T_{6}+3/2T_{1}^{5}T_{2}T_{3}-1/2T_{1}T_{2}T_{3}+3/2T_{1}^{6}T_{3}-1/2T_{1}^{2}T_{3}
-1/2T_{1}^{6}T_{2}+4T_{1}^{4}T_{2}-5/2T_{1}^{2}T_{2}-1/2T_{1}^{7}+2T_{1}^{5}+3/2T_{1}^{3}-2T_{1}$, \\
$\widehat{f_7} = T_{7}+1/2T_{1}^{7}T_{2}T_{3}-2T_{1}^{5}T_{2}T_{3}+1/2T_{1}^{3}T_{2}T_{3}
+4T_{1}^{6}T_{3}+7/2T_{1}^{4}T_{3}-2T_{1}^{2}T_{3}-1/2T_{3}+7/2T_{1}^{4}T_{2}-2T_{1}^{2}T_{2}
-1/2T_{2}-1/2T_{1}^{5}-2T_{1}^{3}-1/2T_{1}$, \\
$\widehat{f_8} = T_{8}+5/4T_{1}^{6}T_{2}T_{3}+3/4T_{1}^{4}T_{2}T_{3}-7/4T_{1}^{2}T_{2}T_{3}
-1/4T_{2}T_{3}+1/4T_{1}^{7}T_{2}+1/4T_{1}^{5}T_{2}-3/4T_{1}^{3}T_{2}-3/4T_{1}T_{2}+5/2T_{1}^{4}+4T_{1}^{2}+1/2$.
\end{example}

\begin{example}
$G = 8T24 = E(8):D_6$, 
$p=601$, \\$x=(578, 499, 355, 408, 383, 392, 156, 235)$, $d(f)=54760000$, \\
$\widehat{f_1} = T_{1}^{8}-T_{1}^{7}+T_{1}^{6}+T_{1}^{2}+T_{1}+1$, \\
$\widehat{f_2} = T_{2}^{3}+3/74T_{1}^{7}T_{2}^{2}-25/37T_{1}^{6}T_{2}^{2}+17/37T_{1}^{5}T_{2}^{2}+17/37T_{1}^{
3}T_{2}^{2}+8/37T_{1}^{2}T_{2}^{2}-31/74T_{1}T_{2}^{2}-17/37T_{2}^{2}-23/74T_{1}^{7}T_{2}+14/
37T_{1}^{6}T_{2}+9/74T_{1}^{5}T_{2}-15/74T_{1}^{4}T_{2}-7/74T_{1}^{3}T_{2}-31/74T_{1}^{2}T_{
2}-20/37T_{1}T_{2}+41/74T_{2}+27/74T_{1}^{7}-13/74T_{1}^{6}+67/74T_{1}^{5}+9/74T_{1}^{4}-79/
74T_{1}^{3}-29/74T_{1}^{2}+14/37T_{1}+3/37$, \\
$\widehat{f_3} = T_{3}^{2}-47/74T_{1}^{7}T_{2}^{2}T_{3}+45/37T_{1}^{6}T_{2}^{2}T_{3}-25/74T_{1}^{5}T_{2}^{
2}T_{3}+47/74T_{1}^{4}T_{2}^{2}T_{3}
-1/74T_{1}^{3}T_{2}^{2}T_{3}+43/74T_{1}^{2}T_{2}^{2}T_{
3}+29/37T_{1}T_{2}^{2}T_{3}+77/74T_{2}^{2}T_{3}+58/37T_{1}^{7}T_{2}T_{3}-10/37T_{1}^{6}T_{
2}T_{3}+34/37T_{1}^{5}T_{2}T_{3}-22/37T_{1}^{4}T_{2}T_{3}-68/37T_{1}^{3}T_{2}T_{3}+31/37T_{
1}^{2}T_{2}T_{3}+46/37T_{1}T_{2}T_{3}+28/37T_{2}T_{3}+19/74T_{1}^{7}T_{3}-13/37T_{1}^{6}T_{
3}+29/74T_{1}^{5}T_{3}-149/74T_{1}^{4}T_{3}+39/74T_{1}^{3}T_{3}+77/74T_{1}^{2}T_{3}-15/37T_{
1}T_{3}+17/74T_{3}+31/74T_{1}^{7}T_{2}^{2}+25/37T_{1}^{6}T_{2}^{2}+41/74T_{1}^{5}T_{2}^{2}
-31/74T_{1}^{4}T_{2}^{2}-139/74T_{1}^{3}T_{2}^{2}-19/74T_{1}^{2}T_{2}^{2}+18/37T_{1}T_{2}^{
2}+31/74T_{2}^{2}-71/74T_{1}^{7}T_{2}+53/37T_{1}^{6}T_{2}-23/74T_{1}^{5}T_{2}-175/74T_{1}^{
4}T_{2}+53/74T_{1}^{3}T_{2}+109/74T_{1}^{2}T_{2}+55/74T_{2}+39/37T_{1}^{7}-8/37T_{1}^{6}+13/
37T_{1}^{5}+32/37T_{1}^{4}-12/37T_{1}^{3}-19/37T_{1}^{2}+55/37T_{1}+27/37$, \\
$\widehat{f_4} = T_{4}-9/5T_{1}^{7}T_{2}^{2}T_{3}+3T_{1}^{6}T_{2}^{2}T_{3}+3/5T_{1}^{5}T_{2}^{2}T_{3}+37/10T_{
1}^{4}T_{2}^{2}T_{3}-43/10T_{1}^{3}T_{2}^{2}T_{3}-8/5T_{1}^{2}T_{2}^{2}T_{3}+1/10T_{1}T_{2}^{
2}T_{3}+7/10T_{2}^{2}T_{3}-3/5T_{1}^{7}T_{2}T_{3}+16/5T_{1}^{6}T_{2}T_{3}-1/2T_{1}^{5}T_{
2}T_{3}-27/5T_{1}^{4}T_{2}T_{3}-38/5T_{1}^{3}T_{2}T_{3}+41/10T_{1}^{2}T_{2}T_{3}-1/10T_{1}T_{
2}T_{3}+5/2T_{2}T_{3}+7/5T_{1}^{7}T_{3}+1/5T_{1}^{6}T_{3}-9/10T_{1}^{5}T_{3}-33/10T_{1}^{
4}T_{3}+33/10T_{1}^{3}T_{3}+63/10T_{1}^{2}T_{3}+21/5T_{1}T_{3}+26/5T_{3}-16/5T_{1}^{7}T_{2}^{
2}+17/5T_{1}^{6}T_{2}^{2}-9/10T_{1}^{5}T_{2}^{2}-9/10T_{1}^{4}T_{2}^{2}-13/2T_{1}^{3}T_{2}^{
2}-19/10T_{1}^{2}T_{2}^{2}-17/5T_{1}T_{2}^{2}-7/5T_{2}^{2}-9/10T_{1}^{7}T_{2}+43/10T_{1}^{
6}T_{2}-29/10T_{1}^{5}T_{2}-12/5T_{1}^{4}T_{2}+1/5T_{1}^{3}T_{2}+79/10T_{1}^{2}T_{2}+28/5T_{
1}T_{2}+43/5T_{2}+51/10T_{1}^{7}-29/10T_{1}^{6}+16/5T_{1}^{5}-1/2T_{1}^{4}+29/10T_{1}^{3}+19/
5T_{1}^{2}+5T_{1}+13/5$, \\
$\widehat{f_5} = T_{5}+7/5T_{1}^{7}T_{2}^{2}T_{3}-13/5T_{1}^{6}T_{2}^{2}T_{3}-2/5T_{1}^{5}T_{2}^{2}T_{3}+16/
5T_{1}^{4}T_{2}^{2}T_{3}+22/5T_{1}^{3}T_{2}^{2}T_{3}+T_{1}^{2}T_{2}^{2}T_{3}+8/5T_{1}T_{2}^{
2}T_{3}-9/5T_{2}^{2}T_{3}+6/5T_{1}^{7}T_{2}T_{3}-18/5T_{1}^{6}T_{2}T_{3}-3/5T_{1}^{5}T_{2}T_{
3}+6/5T_{1}^{4}T_{2}T_{3}+26/5T_{1}^{3}T_{2}T_{3}+2/5T_{1}^{2}T_{2}T_{3}-12/5T_{1}T_{2}T_{3}
-11/5T_{2}T_{3}+T_{1}^{7}T_{3}-11/5T_{1}^{6}T_{3}+4/5T_{1}^{5}T_{3}+4/5T_{1}^{4}T_{3}+16/5T_{
1}^{3}T_{3}+1/5T_{1}T_{3}-3/5T_{3}+T_{1}^{7}T_{2}^{2}-9/5T_{1}^{6}T_{2}^{2}-2/5T_{1}^{5}T_{
2}^{2}+6/5T_{1}^{4}T_{2}^{2}+6/5T_{1}^{3}T_{2}^{2}-14/5T_{1}^{2}T_{2}^{2}-17/5T_{1}T_{2}^{2}
-19/5T_{2}^{2}+3T_{1}^{7}T_{2}-24/5T_{1}^{6}T_{2}+2T_{1}^{5}T_{2}+6/5T_{1}^{4}T_{2}+2T_{1}^{
3}T_{2}-11/5T_{1}^{2}T_{2}-12/5T_{1}T_{2}-12/5T_{2}+7/5T_{1}^{7}-8/5T_{1}^{6}+3/5T_{1}^{5}+4/
5T_{1}^{4}+6/5T_{1}^{3}-8/5T_{1}^{2}-T_{1}-3/5$, \\
$\widehat{f_6} = T_{6}-1/5T_{1}^{7}T_{2}^{2}T_{3}-9/5T_{1}^{6}T_{2}^{2}T_{3}-14/5T_{1}^{5}T_{2}^{2}T_{3}-3/
5T_{1}^{4}T_{2}^{2}T_{3}-3/5T_{1}^{3}T_{2}^{2}T_{3}-2/5T_{1}^{2}T_{2}^{2}T_{3}-2/5T_{1}T_{
2}^{2}T_{3}-2/5T_{2}^{2}T_{3}+4/5T_{1}^{7}T_{2}T_{3}-12/5T_{1}^{6}T_{2}T_{3}+4/5T_{1}^{4}T_{
2}T_{3}+8/5T_{1}^{3}T_{2}T_{3}-2/5T_{1}^{2}T_{2}T_{3}-8/5T_{1}T_{2}T_{3}-2T_{2}T_{3}-2T_{1}^{
7}T_{3}-4/5T_{1}^{6}T_{3}-3T_{1}^{5}T_{3}+1/5T_{1}^{4}T_{3}+2/5T_{1}^{3}T_{3}+1/5T_{1}^{2}T_{
3}+1/5T_{1}T_{3}+12/5T_{1}^{7}T_{2}^{2}-T_{1}^{6}T_{2}^{2}+1/5T_{1}^{5}T_{2}^{2}+2/5T_{1}^{
4}T_{2}^{2}+1/5T_{1}^{3}T_{2}^{2}-3T_{1}^{2}T_{2}^{2}-16/5T_{1}T_{2}^{2}-16/5T_{2}^{2}+2/5T_{
1}^{7}T_{2}-2T_{1}^{6}T_{2}-6/5T_{1}^{5}T_{2}+4/5T_{1}^{4}T_{2}+4/5T_{1}^{3}T_{2}-4/5T_{1}^{
2}T_{2}-4/5T_{1}T_{2}-2/5T_{2}+4/5T_{1}^{7}-8/5T_{1}^{6}+2/5T_{1}^{5}-3/5T_{1}^{4}-3/5T_{1}^{
3}-14/5T_{1}^{2}-11/5T_{1}-11/5$, \\
$\widehat{f_7} = T_{7}+13/5T_{1}^{7}T_{2}^{2}T_{3}-21/10T_{1}^{6}T_{2}^{2}T_{3}-13/20T_{1}^{5}T_{2}^{2}T_{3}
-3/2T_{1}^{4}T_{2}^{2}T_{3}+49/10T_{1}^{3}T_{2}^{2}T_{3}+51/20T_{1}^{2}T_{2}^{2}T_{3}+7/4T_{
1}T_{2}^{2}T_{3}+21/20T_{2}^{2}T_{3}+9/4T_{1}^{7}T_{2}T_{3}-117/20T_{1}^{6}T_{2}T_{3}+8/5T_{
1}^{5}T_{2}T_{3}+23/10T_{1}^{4}T_{2}T_{3}+31/10T_{1}^{3}T_{2}T_{3}-27/5T_{1}^{2}T_{2}T_{3}
-51/20T_{1}T_{2}T_{3}-109/20T_{2}T_{3}-67/20T_{1}^{7}T_{3}+7/20T_{1}^{6}T_{3}-27/20T_{1}^{
5}T_{3}+2T_{1}^{4}T_{3}-2/5T_{1}^{3}T_{3}-91/20T_{1}^{2}T_{3}-9/2T_{1}T_{3}-24/5T_{3}-16/5T_{
1}^{7}T_{2}^{2}-3/10T_{1}^{6}T_{2}^{2}-47/20T_{1}^{5}T_{2}^{2}+18/5T_{1}^{4}T_{2}^{2}+8/5T_{
1}^{3}T_{2}^{2}-113/20T_{1}^{2}T_{2}^{2}-69/20T_{1}T_{2}^{2}-89/20T_{2}^{2}-127/20T_{1}^{
7}T_{2}+113/20T_{1}^{6}T_{2}-T_{1}^{5}T_{2}+39/10T_{1}^{4}T_{2}-63/10T_{1}^{3}T_{2}-12/5T_{
1}^{2}T_{2}-39/20T_{1}T_{2}+1/20T_{2}-71/20T_{1}^{7}+121/20T_{1}^{6}-33/20T_{1}^{5}-1/10T_{
1}^{4}-51/10T_{1}^{3}+33/20T_{1}^{2}+7/10T_{1}+11/5$, \\
$\widehat{f_8} = T_{8}+29/74T_{1}^{7}T_{2}^{2}T_{3}+12/37T_{1}^{6}T_{2}^{2}T_{3}+97/74T_{1}^{5}T_{2}^{2}T_{3}
-127/74T_{1}^{4}T_{2}^{2}T_{3}-139/74T_{1}^{3}T_{2}^{2}T_{3}-69/74T_{1}^{2}T_{2}^{2}T_{3}-39/
37T_{1}T_{2}^{2}T_{3}-93/74T_{2}^{2}T_{3}-91/37T_{1}^{7}T_{2}T_{3}+83/37T_{1}^{6}T_{2}T_{3}
-60/37T_{1}^{5}T_{2}T_{3}-18/37T_{1}^{4}T_{2}T_{3}+6/37T_{1}^{3}T_{2}T_{3}-26/37T_{1}^{2}T_{
2}T_{3}+17/37T_{1}T_{2}T_{3}+43/37T_{2}T_{3}-109/74T_{1}^{7}T_{3}+79/37T_{1}^{6}T_{3}-81/
74T_{1}^{5}T_{3}+123/74T_{1}^{4}T_{3}-25/74T_{1}^{3}T_{3}-45/74T_{1}^{2}T_{3}+30/37T_{1}T_{
3}+21/74T_{3}-27/37T_{1}^{7}T_{2}^{2}+60/37T_{1}^{6}T_{2}^{2}+7/37T_{1}^{5}T_{2}^{2}-9/37T_{
1}^{4}T_{2}^{2}-71/37T_{1}^{3}T_{2}^{2}-50/37T_{1}^{2}T_{2}^{2}-10/37T_{1}T_{2}^{2}+3/37T_{
2}^{2}-72/37T_{1}^{7}T_{2}+121/37T_{1}^{6}T_{2}-7/37T_{1}^{5}T_{2}-25/37T_{1}^{4}T_{2}-25/
37T_{1}^{3}T_{2}+29/37T_{1}^{2}T_{2}+23/37T_{1}T_{2}+58/37T_{2}+95/74T_{1}^{7}+7/37T_{1}^{6}+
103/74T_{1}^{5}+25/74T_{1}^{4}-173/74T_{1}^{3}-67/74T_{1}^{2}+42/37T_{1}+35/74$.
\end{example}

\begin{example}
$G = 9T11 = E(9):6=1/2[3^2:2]S(3)$, 
$p=991$, \\$x=(981, 149, 630, 523, 829, 772, 191, 28, 852)$, $d(f)=13177032454057536$, \\
$\widehat{f_1} = T_{1}^{9}-3T_{1}^{6}+3T_{1}^{3}+8$, \\
$\widehat{f_2} = T_{2}^{2}+T_{1}^{2}T_{2}-T_{1}^{3}$, \\
$\widehat{f_3} = T_{3}^{3}+1/9T_{1}^{7}T_{2}T_{3}^{2}-2/9T_{1}^{4}T_{2}T_{3}^{2}-8/9T_{1}T_{2}T_{3}^{2}+2/
9T_{1}^{8}T_{3}^{2}-4/9T_{1}^{5}T_{3}^{2}-16/9T_{1}^{2}T_{3}^{2}-2/9T_{1}^{8}T_{2}T_{3}+4/
9T_{1}^{5}T_{2}T_{3}+7/9T_{1}^{2}T_{2}T_{3}-1/9T_{1}^{6}T_{3}+20/9T_{1}^{3}T_{3}+8/9T_{3}+
1/9T_{1}^{6}T_{2}-11/9T_{1}^{3}T_{2}-8/9T_{2}-1/9T_{1}^{7}-16/9T_{1}^{4}+8/9T_{1}$, \\
$\widehat{f_4} = T_{4}+1/3T_{1}^{7}T_{3}^{2}-2/3T_{1}^{4}T_{3}^{2}-8/3T_{1}T_{3}^{2}+1/3T_{1}^{7}T_{2}T_{3}
-2/3T_{1}^{4}T_{2}T_{3}-8/3T_{1}T_{2}T_{3}-1/3T_{1}^{8}T_{2}+2/3T_{1}^{5}T_{2}-1/3T_{1}^{
2}T_{2}-1/3T_{1}^{6}+11/3T_{1}^{3}+8/3$, \\
$\widehat{f_5} = T_{5}-1/3T_{1}^{6}T_{2}T_{3}^{2}+2/3T_{1}^{3}T_{2}T_{3}^{2}+8/3T_{2}T_{3}^{2}-1/3T_{1}^{
7}T_{2}T_{3}+2/3T_{1}^{4}T_{2}T_{3}+8/3T_{1}T_{2}T_{3}-1/3T_{1}^{8}T_{2}+2/3T_{1}^{5}T_{
2}+8/3T_{1}^{2}T_{2}+3T_{1}^{3}$, \\
$\widehat{f_6} = T_{6}+T_{1}^{5}T_{2}T_{3}^{2}-T_{1}^{2}T_{2}T_{3}^{2}-1/6T_{1}^{6}T_{2}T_{3}+1/3T_{1}^{
3}T_{2}T_{3}+4/3T_{2}T_{3}-2/3T_{1}^{7}T_{3}+5/6T_{1}^{4}T_{3}+4/3T_{1}T_{3}-T_{1}^{4}T_{
2}+4T_{1}T_{2}+1/3T_{1}^{8}-5/3T_{1}^{5}+4/3T_{1}^{2}$, \\
$\widehat{f_7} = T_{7}-1/3T_{1}^{5}T_{2}T_{3}^{2}+1/3T_{1}^{2}T_{2}T_{3}^{2}-2/3T_{1}^{6}T_{3}^{2}+2/3T_{
1}^{3}T_{3}^{2}-2/3T_{1}^{6}T_{2}T_{3}+1/2T_{1}^{3}T_{2}T_{3}-4/3T_{2}T_{3}-1/6T_{1}^{
7}T_{3}+1/6T_{1}^{4}T_{3}+5/3T_{1}^{4}T_{2}+4/3T_{1}T_{2}+1/3T_{1}^{8}-1/3T_{1}^{5}$, \\
$\widehat{f_8} = T_{8}+1/3T_{1}^{5}T_{2}T_{3}^{2}-1/3T_{1}^{2}T_{2}T_{3}^{2}+2/3T_{1}^{6}T_{3}^{2}-2/3T_{
1}^{3}T_{3}^{2}+5/6T_{1}^{6}T_{2}T_{3}-5/6T_{1}^{3}T_{2}T_{3}+5/6T_{1}^{7}T_{3}-T_{1}^{
4}T_{3}-4/3T_{1}T_{3}+1/3T_{1}^{4}T_{2}+8/3T_{1}T_{2}+1/3T_{1}^{8}+T_{1}^{5}-4/3T_{1}^{2}$, \\
$\widehat{f_9} = T_{9}-2/3T_{1}^{6}T_{2}T_{3}^{2}+4/3T_{1}^{3}T_{2}T_{3}^{2}+16/3T_{2}T_{3}^{2}-1/3T_{1}^{
7}T_{3}^{2}+2/3T_{1}^{4}T_{3}^{2}+8/3T_{1}T_{3}^{2}-1/3T_{1}^{8}T_{2}+2/3T_{1}^{5}T_{2}
-19/3T_{1}^{2}T_{2}-2/3T_{1}^{6}-5/3T_{1}^{3}+16/3$.
\end{example}

\begin{example}
$G = 9T18 = E(9) : D12 = [32 : 2]S(3) = [1/2.S(3)2]S(3)$, 
$p=1231$, $x=(1195, 879, 893, 1072, 497, 391, 814, 26, 388)$, $d(f)=239483061$, \\
$\widehat{f_1} = T_{1}^{9}-T_{1}^{3}-1$, \\
$\widehat{f_2} = T_{2}^{2}+T_{1}^{2}T_{2}-T_{1}^{3}$, \\
$\widehat{f_3} = T_{3}^{6}-6/23T_{1}^{8}T_{3}^{5}-45/23T_{1}^{5}T_{3}^{5}-54/23T_{1}^{2}T_{3}^{5}-6/23T_{1}^{
8}T_{2}T_{3}^{4}+9/23T_{1}^{5}T_{2}T_{3}^{4}+4/23T_{1}^{2}T_{2}T_{3}^{4}+27/23T_{1}^{6}T_{
3}^{4}+29/23T_{1}^{3}T_{3}^{4}+9/23T_{1}^{6}T_{2}T_{3}^{3}+33/23T_{1}^{3}T_{2}T_{3}^{3}+12/
23T_{2}T_{3}^{3}-36/23T_{1}^{8}T_{3}^{2}-27/23T_{1}^{5}T_{3}^{2}+6/23T_{1}^{2}T_{3}^{2}+18/
23T_{1}^{8}T_{2}T_{3}-27/23T_{1}^{5}T_{2}T_{3}-12/23T_{1}^{2}T_{2}T_{3}+18/23T_{1}^{3}T_{3}+
27/23T_{3}+54/23T_{1}^{6}T_{2}+6/23T_{1}^{3}T_{2}-9/23T_{2}$, \\
$\widehat{f_4} = T_{4}-27/23T_{1}^{6}T_{2}T_{3}^{5}-30/23T_{1}^{3}T_{2}T_{3}^{5}-36/23T_{2}T_{3}^{5}+27/23T_{
1}^{7}T_{3}^{5}+30/23T_{1}^{4}T_{3}^{5}+36/23T_{1}T_{3}^{5}+54/23T_{1}^{7}T_{2}T_{3}^{4}+60/
23T_{1}^{4}T_{2}T_{3}^{4}+72/23T_{1}T_{2}T_{3}^{4}+27/23T_{1}^{8}T_{3}^{4}+30/23T_{1}^{5}T_{
3}^{4}+36/23T_{1}^{2}T_{3}^{4}+54/23T_{1}^{8}T_{2}T_{3}^{3}+6/23T_{1}^{5}T_{2}T_{3}^{3}-9/
23T_{1}^{2}T_{2}T_{3}^{3}-24/23T_{1}^{6}T_{3}^{3}-18/23T_{1}^{3}T_{3}^{3}+27/23T_{3}^{3}-24/
23T_{1}^{6}T_{2}T_{3}^{2}-18/23T_{1}^{3}T_{2}T_{3}^{2}+27/23T_{2}T_{3}^{2}+24/23T_{1}^{7}T_{
3}^{2}+18/23T_{1}^{4}T_{3}^{2}-27/23T_{1}T_{3}^{2}+48/23T_{1}^{7}T_{2}T_{3}+36/23T_{1}^{4}T_{
2}T_{3}-54/23T_{1}T_{2}T_{3}+24/23T_{1}^{8}T_{3}+18/23T_{1}^{5}T_{3}-27/23T_{1}^{2}T_{3}-6/
23T_{1}^{8}T_{2}+36/23T_{1}^{5}T_{2}+33/23T_{1}^{2}T_{2}+18/23T_{1}^{6}+30/23T_{1}^{3}-30/23$, \\
$\widehat{f_5} = T_{5}-54/23T_{1}^{6}T_{2}T_{3}^{5}-60/23T_{1}^{3}T_{2}T_{3}^{5}-72/23T_{2}T_{3}^{5}-27/23T_{
1}^{7}T_{3}^{5}-30/23T_{1}^{4}T_{3}^{5}-36/23T_{1}T_{3}^{5}-54/23T_{1}^{7}T_{2}T_{3}^{4}-60/
23T_{1}^{4}T_{2}T_{3}^{4}-72/23T_{1}T_{2}T_{3}^{4}-27/23T_{1}^{8}T_{3}^{4}-30/23T_{1}^{5}T_{
3}^{4}-36/23T_{1}^{2}T_{3}^{4}+27/23T_{1}^{8}T_{2}T_{3}^{3}-24/23T_{1}^{5}T_{2}T_{3}^{3}-45/
23T_{1}^{2}T_{2}T_{3}^{3}-30/23T_{1}^{6}T_{3}^{3}-63/23T_{1}^{3}T_{3}^{3}-27/23T_{3}^{3}-48/
23T_{1}^{6}T_{2}T_{3}^{2}-36/23T_{1}^{3}T_{2}T_{3}^{2}+54/23T_{2}T_{3}^{2}-24/23T_{1}^{7}T_{
3}^{2}-18/23T_{1}^{4}T_{3}^{2}+27/23T_{1}T_{3}^{2}-48/23T_{1}^{7}T_{2}T_{3}-36/23T_{1}^{4}T_{
2}T_{3}+54/23T_{1}T_{2}T_{3}-24/23T_{1}^{8}T_{3}-18/23T_{1}^{5}T_{3}+27/23T_{1}^{2}T_{3}-30/
23T_{1}^{8}T_{2}+18/23T_{1}^{5}T_{2}+60/23T_{1}^{2}T_{2}-18/23T_{1}^{6}+3/23T_{1}^{3}-24/23$, \\
$\widehat{f_6} = T_{6}-45/23T_{1}^{6}T_{2}T_{3}^{5}-27/23T_{1}^{3}T_{2}T_{3}^{5}-60/23T_{2}T_{3}^{5}-18/23T_{
1}^{7}T_{2}T_{3}^{4}-8/23T_{1}^{4}T_{2}T_{3}^{4}-6/23T_{1}T_{2}T_{3}^{4}-42/23T_{1}^{8}T_{
3}^{4}-37/23T_{1}^{5}T_{3}^{4}+39/23T_{1}^{2}T_{3}^{4}-32/23T_{1}^{8}T_{2}T_{3}^{3}-6/23T_{
1}^{5}T_{2}T_{3}^{3}+63/23T_{1}^{2}T_{2}T_{3}^{3}-18/23T_{1}^{6}T_{3}^{3}+41/23T_{1}^{3}T_{
3}^{3}-13/23T_{3}^{3}-33/23T_{1}^{6}T_{2}T_{3}^{2}-75/23T_{1}^{3}T_{2}T_{3}^{2}-90/23T_{2}T_{
3}^{2}-27/23T_{1}^{7}T_{2}T_{3}-12/23T_{1}^{4}T_{2}T_{3}-9/23T_{1}T_{2}T_{3}-63/23T_{1}^{
8}T_{3}-21/23T_{1}^{5}T_{3}+24/23T_{1}^{2}T_{3}-48/23T_{1}^{8}T_{2}-9/23T_{1}^{5}T_{2}+37/
23T_{1}^{2}T_{2}-27/23T_{1}^{6}-19/23T_{1}^{3}-54/23$, \\
$\widehat{f_7} = T_{7}-9/23T_{1}^{6}T_{2}T_{3}^{5}-33/23T_{1}^{3}T_{2}T_{3}^{5}-12/23T_{2}T_{3}^{5}+18/23T_{
1}^{7}T_{3}^{5}-3/23T_{1}^{4}T_{3}^{5}+24/23T_{1}T_{3}^{5}+24/23T_{1}^{7}T_{2}T_{3}^{4}+7/
23T_{1}^{4}T_{2}T_{3}^{4}-9/23T_{1}T_{2}T_{3}^{4}+42/23T_{1}^{8}T_{3}^{4}+26/23T_{1}^{5}T_{
3}^{4}-21/23T_{1}^{2}T_{3}^{4}+16/23T_{1}^{8}T_{2}T_{3}^{3}+30/23T_{1}^{5}T_{2}T_{3}^{3}+9/
23T_{1}^{2}T_{2}T_{3}^{3}+38/23T_{1}^{3}T_{3}^{3}+11/23T_{3}^{3}+21/23T_{1}^{6}T_{2}T_{3}^{2}
-15/23T_{1}^{3}T_{2}T_{3}^{2}-18/23T_{2}T_{3}^{2}+27/23T_{1}^{7}T_{3}^{2}+30/23T_{1}^{4}T_{
3}^{2}+36/23T_{1}T_{3}^{2}+36/23T_{1}^{7}T_{2}T_{3}+45/23T_{1}^{4}T_{2}T_{3}+21/23T_{1}T_{
2}T_{3}+63/23T_{1}^{8}T_{3}+39/23T_{1}^{5}T_{3}+3/23T_{1}^{2}T_{3}+24/23T_{1}^{8}T_{2}+45/
23T_{1}^{5}T_{2}+25/23T_{1}^{2}T_{2}+11/23T_{1}^{3}-18/23$, \\
$\widehat{f_8} = T_{8}-27/23T_{1}^{6}T_{2}T_{3}^{5}-30/23T_{1}^{3}T_{2}T_{3}^{5}-36/23T_{2}T_{3}^{5}-18/23T_{
1}^{7}T_{3}^{5}+3/23T_{1}^{4}T_{3}^{5}-24/23T_{1}T_{3}^{5}-6/23T_{1}^{7}T_{2}T_{3}^{4}+1/
23T_{1}^{4}T_{2}T_{3}^{4}+15/23T_{1}T_{2}T_{3}^{4}+11/23T_{1}^{5}T_{3}^{4}-18/23T_{1}^{2}T_{
3}^{4}-8/23T_{1}^{8}T_{2}T_{3}^{3}+12/23T_{1}^{5}T_{2}T_{3}^{3}+36/23T_{1}^{2}T_{2}T_{3}^{3}+
18/23T_{1}^{6}T_{3}^{3}+35/23T_{1}^{3}T_{3}^{3}+35/23T_{3}^{3}-6/23T_{1}^{6}T_{2}T_{3}^{2}
-45/23T_{1}^{3}T_{2}T_{3}^{2}-54/23T_{2}T_{3}^{2}-27/23T_{1}^{7}T_{3}^{2}-30/23T_{1}^{4}T_{
3}^{2}-36/23T_{1}T_{3}^{2}-9/23T_{1}^{7}T_{2}T_{3}-33/23T_{1}^{4}T_{2}T_{3}-12/23T_{1}T_{
2}T_{3}-18/23T_{1}^{5}T_{3}-27/23T_{1}^{2}T_{3}-12/23T_{1}^{8}T_{2}+18/23T_{1}^{5}T_{2}+31/
23T_{1}^{2}T_{2}+27/23T_{1}^{6}+41/23T_{1}^{3}+18/23$, \\
$\widehat{f_9} = T_{9}+162/23T_{1}^{6}T_{2}T_{3}^{5}-36/23T_{1}^{3}T_{2}T_{3}^{5}-108/23T_{2}T_{3}^{5}-81/
23T_{1}^{8}T_{2}T_{3}^{3}+18/23T_{1}^{5}T_{2}T_{3}^{3}+54/23T_{1}^{2}T_{2}T_{3}^{3}+18/23T_{
1}^{6}T_{3}^{3}-27/23T_{1}^{3}T_{3}^{3}-81/23T_{3}^{3}-18/23T_{1}^{6}T_{2}T_{3}^{2}+27/23T_{
1}^{3}T_{2}T_{3}^{2}+81/23T_{2}T_{3}^{2}-72/23T_{1}^{8}T_{2}+108/23T_{1}^{5}T_{2}-90/23T_{
1}^{2}T_{2}+108/23T_{1}^{6}-162/23T_{1}^{3}-72/23$.
\end{example}

\begin{example}
$G = 10T5 = F(5)[x]2$, 
$p=601$, \\$x=(572, 533, 354, 284, 93, 228, 274, 73, 510, 84)$, $d(f)=320000000000$, \\
$\widehat{f_1} = T_{1}^{10}-2T_{1}^{5}-1$, \\
$\widehat{f_2} = T_{2}^{4}-T_{1}^{4}T_{2}+T_{1}^{5}$, \\
$\widehat{f_3} = T_{3}+5/4T_{1}^{9}T_{2}^{3}-5/4T_{1}^{7}T_{2}^{3}+5/4T_{1}^{5}T_{2}^{3}-5/4T_{1}^{4}T_{
2}^{3}+5/4T_{1}^{2}T_{2}^{3}-5/4T_{2}^{3}+5/4T_{1}^{9}T_{2}^{2}-5/4T_{1}^{8}T_{2}^{2}+5/
2T_{1}^{6}T_{2}^{2}-15/4T_{1}^{4}T_{2}^{2}+5/4T_{1}^{3}T_{2}^{2}-5/2T_{1}T_{2}^{2}-5/4T_{
1}^{8}T_{2}+5/4T_{1}^{7}T_{2}-5/4T_{1}^{5}T_{2}+15/4T_{1}^{3}T_{2}-5/4T_{1}^{2}T_{2}+5/
4T_{2}$, \\
$\widehat{f_4} = T_{4}-T_{1}^{8}T_{2}^{3}+4T_{1}^{5}T_{2}^{3}+2T_{2}^{3}-2T_{1}^{9}T_{2}^{2}+T_{1}^{7}T_{
2}^{2}+T_{1}^{8}T_{2}+T_{1}^{6}T_{2}-T_{1}^{5}$, \\
$\widehat{f_5} = T_{5}+5/4T_{1}^{9}T_{2}^{3}+5/4T_{1}^{7}T_{2}^{3}-5/4T_{1}^{4}T_{2}^{3}-5/4T_{1}^{2}T_{
2}^{3}+5/4T_{1}^{9}T_{2}^{2}-5/4T_{1}^{6}T_{2}^{2}-15/4T_{1}^{4}T_{2}^{2}+5/4T_{1}T_{2}^{
2}$, \\
$\widehat{f_6} = T_{6}+T_{1}^{9}T_{2}^{3}+2T_{1}^{6}T_{2}^{3}+T_{1}T_{2}^{3}-T_{1}^{8}T_{2}^{2}-2T_{1}^{
5}T_{2}^{2}-T_{2}^{2}+T_{1}^{9}T_{2}+T_{1}^{7}T_{2}-T_{1}^{8}-T_{1}^{6}$, \\
$\widehat{f_7} = T_{7}+5/4T_{1}^{9}T_{2}^{3}-5/4T_{1}^{4}T_{2}^{3}-5/4T_{1}^{9}T_{2}^{2}-5/4T_{1}^{8}T_{
2}^{2}+15/4T_{1}^{4}T_{2}^{2}+5/4T_{1}^{3}T_{2}^{2}+5/4T_{1}^{7}T_{2}-5/4T_{1}^{5}T_{2}-5/
4T_{1}^{2}T_{2}+5/4T_{2}-5/4T_{1}^{9}-5/4T_{1}^{6}+15/4T_{1}^{4}+5/4T_{1}$, \\
$\widehat{f_8} = T_{8}-T_{1}^{8}T_{2}^{3}+2T_{1}^{5}T_{2}^{3}+T_{2}^{3}-T_{1}^{8}T_{2}-T_{1}^{5}$, \\
$\widehat{f_9} = T_{9}+5/4T_{1}^{9}T_{2}^{3}-5/4T_{1}^{5}T_{2}^{3}-5/4T_{1}^{4}T_{2}^{3}+5/4T_{2}^{3}-5/
4T_{1}^{9}T_{2}^{2}-5/4T_{1}^{8}T_{2}^{2}-5/4T_{1}^{6}T_{2}^{2}+15/4T_{1}^{4}T_{2}^{2}+5/
4T_{1}^{3}T_{2}^{2}+5/4T_{1}T_{2}^{2}+5/4T_{1}^{8}T_{2}+5/2T_{1}^{5}T_{2}-15/4T_{1}^{3}T_{
2}-5/2T_{2}+5/4T_{1}^{9}-15/4T_{1}^{4}$, \\
$\widehat{f_{10}} = T_{10}-T_{1}^{8}T_{2}^{3}+2T_{1}^{5}T_{2}^{3}+T_{2}^{3}-T_{1}^{9}T_{2}^{2}+2T_{1}^{8}T_{
2}-T_{1}^{7}-T_{1}^{5}$.
\end{example}

\begin{example}
$G = 10T11 = A(5)[X]2$, 
$p=3253$, \\$x=(3187, 2586, 2354, 2647, 1290, 66, 667, 899, 606, 1963)$, \\$d(f)=-5739519416467456$, \\
$\widehat{f_1} = T_{1}^{10}+T_{1}^{8}-4T_{1}^{2}+4$, \\
$\widehat{f_2} = T_{2}^{4}+35/136T_{1}^{8}T_{2}^{3}-7/34T_{1}^{6}T_{2}^{3}+7/34T_{1}^{4}T_{2}^{3}+7/17T_{1}^{
2}T_{2}^{3}-7/34T_{2}^{3}-1/68T_{1}^{9}T_{2}^{2}+13/34T_{1}^{7}T_{2}^{2}-7/34T_{1}^{5}T_{2}^{
2}-13/17T_{1}^{3}T_{2}^{2}+8/17T_{1}T_{2}^{2}+3/34T_{1}^{6}T_{2}+16/17T_{1}^{4}T_{2}+5/17T_{
1}^{2}T_{2}-4/17T_{2}+11/68T_{1}^{9}+3/34T_{1}^{7}-16/17T_{1}^{5}+11/17T_{1}^{3}-19/17T_{1}$, \\
$\widehat{f_3} = T_{3}^{3}-159/272T_{1}^{9}T_{2}^{3}T_{3}^{2}-6/17T_{1}^{7}T_{2}^{3}T_{3}^{2}+19/272T_{1}^{
5}T_{2}^{3}T_{3}^{2}+31/136T_{1}^{3}T_{2}^{3}T_{3}^{2}-71/136T_{1}T_{2}^{3}T_{3}^{2}-5/17T_{
1}^{8}T_{2}^{2}T_{3}^{2}+235/272T_{1}^{6}T_{2}^{2}T_{3}^{2}-173/136T_{1}^{4}T_{2}^{2}T_{3}^{
2}+331/136T_{1}^{2}T_{2}^{2}T_{3}^{2}-14/17T_{2}^{2}T_{3}^{2}+125/272T_{1}^{9}T_{2}T_{3}^{2}
-107/272T_{1}^{7}T_{2}T_{3}^{2}+7/17T_{1}^{5}T_{2}T_{3}^{2}-177/136T_{1}^{3}T_{2}T_{3}^{2}+
13/34T_{1}T_{2}T_{3}^{2}+141/272T_{1}^{8}T_{3}^{2}-6/17T_{1}^{6}T_{3}^{2}+7/8T_{1}^{4}T_{3}^{
2}-109/68T_{1}^{2}T_{3}^{2}+22/17T_{3}^{2}+27/272T_{1}^{8}T_{2}^{3}T_{3}-23/34T_{1}^{6}T_{
2}^{3}T_{3}+65/136T_{1}^{4}T_{2}^{3}T_{3}+31/68T_{1}^{2}T_{2}^{3}T_{3}-35/34T_{2}^{3}T_{3}-9/
16T_{1}^{9}T_{2}^{2}T_{3}-9/68T_{1}^{7}T_{2}^{2}T_{3}-3/34T_{1}^{5}T_{2}^{2}T_{3}+5/34T_{1}^{
3}T_{2}^{2}T_{3}-15/68T_{1}T_{2}^{2}T_{3}-95/272T_{1}^{8}T_{2}T_{3}+185/136T_{1}^{6}T_{2}T_{
3}-11/17T_{1}^{4}T_{2}T_{3}-43/68T_{1}^{2}T_{2}T_{3}+139/68T_{2}T_{3}+113/136T_{1}^{9}T_{3}+
1/34T_{1}^{7}T_{3}-7/17T_{1}^{5}T_{3}-25/34T_{1}T_{3}-3/17T_{1}^{9}T_{2}^{3}+41/272T_{1}^{
7}T_{2}^{3}-1/4T_{1}^{5}T_{2}^{3}+63/136T_{1}^{3}T_{2}^{3}-39/68T_{1}T_{2}^{3}-3/8T_{1}^{
8}T_{2}^{2}-81/136T_{1}^{6}T_{2}^{2}+31/136T_{1}^{4}T_{2}^{2}+99/34T_{1}^{2}T_{2}^{2}-229/
68T_{2}^{2}-49/68T_{1}^{9}T_{2}-171/136T_{1}^{7}T_{2}+97/136T_{1}^{5}T_{2}-31/34T_{1}^{3}T_{
2}+47/68T_{1}T_{2}+91/136T_{1}^{8}+139/136T_{1}^{6}-49/68T_{1}^{4}-53/68T_{1}^{2}+35/34$, \\
$\widehat{f_4} = T_{4}-15/544T_{1}^{8}T_{2}^{3}T_{3}^{2}+47/272T_{1}^{6}T_{2}^{3}T_{3}^{2}+9/17T_{1}^{4}T_{
2}^{3}T_{3}^{2}-135/136T_{1}^{2}T_{2}^{3}T_{3}^{2}+109/136T_{2}^{3}T_{3}^{2}+13/272T_{1}^{
9}T_{2}^{2}T_{3}^{2}+269/272T_{1}^{7}T_{2}^{2}T_{3}^{2}-15/17T_{1}^{5}T_{2}^{2}T_{3}^{2}-23/
136T_{1}^{3}T_{2}^{2}T_{3}^{2}+1/4T_{1}T_{2}^{2}T_{3}^{2}-69/272T_{1}^{8}T_{2}T_{3}^{2}+2/
17T_{1}^{6}T_{2}T_{3}^{2}-15/136T_{1}^{4}T_{2}T_{3}^{2}+35/34T_{1}^{2}T_{2}T_{3}^{2}-39/34T_{
2}T_{3}^{2}+35/272T_{1}^{9}T_{3}^{2}-87/136T_{1}^{7}T_{3}^{2}+97/136T_{1}^{5}T_{3}^{2}-27/
68T_{1}^{3}T_{3}^{2}+2/17T_{1}T_{3}^{2}+23/272T_{1}^{9}T_{2}^{3}T_{3}-11/68T_{1}^{7}T_{2}^{
3}T_{3}+19/136T_{1}^{5}T_{2}^{3}T_{3}+1/34T_{1}^{3}T_{2}^{3}T_{3}-2/17T_{1}T_{2}^{3}T_{3}+33/
68T_{1}^{8}T_{2}^{2}T_{3}-23/136T_{1}^{6}T_{2}^{2}T_{3}+5/34T_{1}^{4}T_{2}^{2}T_{3}-83/68T_{
1}^{2}T_{2}^{2}T_{3}+41/34T_{2}^{2}T_{3}-9/68T_{1}^{9}T_{2}T_{3}+107/136T_{1}^{7}T_{2}T_{3}
-15/34T_{1}^{5}T_{2}T_{3}+39/68T_{1}^{3}T_{2}T_{3}-7/17T_{1}T_{2}T_{3}-57/68T_{1}^{8}T_{3}
-41/68T_{1}^{6}T_{3}+5/68T_{1}^{4}T_{3}+33/17T_{1}^{2}T_{3}-65/34T_{3}-21/136T_{1}^{8}T_{2}^{
3}+69/136T_{1}^{6}T_{2}^{3}-59/136T_{1}^{4}T_{2}^{3}+27/68T_{1}^{2}T_{2}^{3}-15/68T_{2}^{3}+
55/272T_{1}^{9}T_{2}^{2}-83/68T_{1}^{7}T_{2}^{2}+157/136T_{1}^{5}T_{2}^{2}-31/34T_{1}^{3}T_{
2}^{2}+1/2T_{1}T_{2}^{2}+279/272T_{1}^{8}T_{2}+5/136T_{1}^{6}T_{2}-9/68T_{1}^{4}T_{2}-135/
68T_{1}^{2}T_{2}+149/68T_{2}+7/136T_{1}^{9}+115/136T_{1}^{7}-79/68T_{1}^{5}+71/68T_{1}^{3}-5/
17T_{1}$, \\
$\widehat{f_5} = T_{5}+13/136T_{1}^{8}T_{2}^{3}T_{3}^{2}+11/136T_{1}^{6}T_{2}^{3}T_{3}^{2}-113/136T_{1}^{4}T_{
2}^{3}T_{3}^{2}-123/68T_{1}^{2}T_{2}^{3}T_{3}^{2}+77/34T_{2}^{3}T_{3}^{2}+1/4T_{1}^{9}T_{2}^{
2}T_{3}^{2}+199/272T_{1}^{7}T_{2}^{2}T_{3}^{2}+9/17T_{1}^{5}T_{2}^{2}T_{3}^{2}+67/136T_{1}^{
3}T_{2}^{2}T_{3}^{2}-3/17T_{1}T_{2}^{2}T_{3}^{2}-83/136T_{1}^{8}T_{2}T_{3}^{2}-61/68T_{1}^{
6}T_{2}T_{3}^{2}+73/68T_{1}^{4}T_{2}T_{3}^{2}+45/68T_{1}^{2}T_{2}T_{3}^{2}-29/17T_{2}T_{3}^{
2}+1/34T_{1}^{9}T_{3}^{2}+43/136T_{1}^{7}T_{3}^{2}-55/136T_{1}^{5}T_{3}^{2}-77/68T_{1}^{3}T_{
3}^{2}+95/68T_{1}T_{3}^{2}-19/272T_{1}^{9}T_{2}^{3}T_{3}-19/68T_{1}^{7}T_{2}^{3}T_{3}-173/
136T_{1}^{5}T_{2}^{3}T_{3}-11/17T_{1}^{3}T_{2}^{3}T_{3}+13/17T_{1}T_{2}^{3}T_{3}-101/136T_{
1}^{8}T_{2}^{2}T_{3}+305/136T_{1}^{6}T_{2}^{2}T_{3}-5/17T_{1}^{4}T_{2}^{2}T_{3}+43/68T_{1}^{
2}T_{2}^{2}T_{3}+15/17T_{2}^{2}T_{3}+31/34T_{1}^{9}T_{2}T_{3}+31/136T_{1}^{7}T_{2}T_{3}+12/
17T_{1}^{5}T_{2}T_{3}+93/68T_{1}^{3}T_{2}T_{3}-28/17T_{1}T_{2}T_{3}+63/136T_{1}^{8}T_{3}-105/
68T_{1}^{6}T_{3}+1/4T_{1}^{4}T_{3}-22/17T_{1}^{2}T_{3}+7/17T_{3}+287/272T_{1}^{8}T_{2}^{3}
-129/136T_{1}^{6}T_{2}^{3}+11/34T_{1}^{4}T_{2}^{3}-43/34T_{1}^{2}T_{2}^{3}+47/68T_{2}^{3}-11/
8T_{1}^{9}T_{2}^{2}+11/34T_{1}^{7}T_{2}^{2}-101/136T_{1}^{5}T_{2}^{2}+43/17T_{1}^{3}T_{2}^{2}
-71/68T_{1}T_{2}^{2}-93/136T_{1}^{8}T_{2}+125/68T_{1}^{6}T_{2}-71/68T_{1}^{4}T_{2}+131/34T_{
1}^{2}T_{2}-91/34T_{2}+43/68T_{1}^{9}-121/136T_{1}^{7}+25/68T_{1}^{5}-117/68T_{1}^{3}-3/34T_{
1}$, \\
$\widehat{f_6} = T_{6}-129/136T_{1}^{9}T_{2}^{3}T_{3}^{2}-257/136T_{1}^{7}T_{2}^{3}T_{3}^{2}-33/68T_{1}^{5}T_{
2}^{3}T_{3}^{2}+25/68T_{1}^{3}T_{2}^{3}T_{3}^{2}+123/34T_{1}T_{2}^{3}T_{3}^{2}+26/17T_{1}^{
8}T_{2}^{2}T_{3}^{2}+183/68T_{1}^{6}T_{2}^{2}T_{3}^{2}+197/68T_{1}^{4}T_{2}^{2}T_{3}^{2}+30/
17T_{1}^{2}T_{2}^{2}T_{3}^{2}-139/34T_{2}^{2}T_{3}^{2}-7/17T_{1}^{9}T_{2}T_{3}^{2}-9/17T_{
1}^{7}T_{2}T_{3}^{2}-27/17T_{1}^{5}T_{2}T_{3}^{2}-48/17T_{1}^{3}T_{2}T_{3}^{2}-23/17T_{1}T_{
2}T_{3}^{2}+3/17T_{1}^{8}T_{3}^{2}+15/17T_{1}^{6}T_{3}^{2}+21/17T_{1}^{4}T_{3}^{2}+T_{1}^{
2}T_{3}^{2}-8/17T_{3}^{2}-9/17T_{1}^{8}T_{2}^{3}T_{3}-21/17T_{1}^{6}T_{2}^{3}T_{3}-27/17T_{
1}^{4}T_{2}^{3}T_{3}-9/17T_{1}^{2}T_{2}^{3}T_{3}+30/17T_{2}^{3}T_{3}-11/17T_{1}^{9}T_{2}^{
2}T_{3}-21/34T_{1}^{7}T_{2}^{2}T_{3}+73/34T_{1}^{5}T_{2}^{2}T_{3}+63/17T_{1}^{3}T_{2}^{2}T_{
3}+95/17T_{1}T_{2}^{2}T_{3}+107/34T_{1}^{8}T_{2}T_{3}+175/34T_{1}^{6}T_{2}T_{3}+45/17T_{1}^{
4}T_{2}T_{3}-9/17T_{1}^{2}T_{2}T_{3}-216/17T_{2}T_{3}-27/17T_{1}^{9}T_{3}-61/17T_{1}^{7}T_{3}
-78/17T_{1}^{5}T_{3}-66/17T_{1}^{3}T_{3}+82/17T_{1}T_{3}+7/17T_{1}^{9}T_{2}^{3}+39/34T_{1}^{
7}T_{2}^{3}+43/34T_{1}^{5}T_{2}^{3}+9/17T_{1}^{3}T_{2}^{3}-36/17T_{1}T_{2}^{3}-61/34T_{1}^{
8}T_{2}^{2}-93/34T_{1}^{6}T_{2}^{2}-12/17T_{1}^{4}T_{2}^{2}+24/17T_{1}^{2}T_{2}^{2}+88/17T_{
2}^{2}+177/68T_{1}^{9}T_{2}+377/68T_{1}^{7}T_{2}+213/34T_{1}^{5}T_{2}+159/34T_{1}^{3}T_{2}
-103/17T_{1}T_{2}-10/17T_{1}^{8}-3/2T_{1}^{6}-109/34T_{1}^{4}-108/17T_{1}^{2}-9/17$, \\
$\widehat{f_7} = T_{7}-231/272T_{1}^{9}T_{2}^{3}T_{3}^{2}+29/136T_{1}^{7}T_{2}^{3}T_{3}^{2}-107/136T_{1}^{
5}T_{2}^{3}T_{3}^{2}+123/68T_{1}^{3}T_{2}^{3}T_{3}^{2}-19/17T_{1}T_{2}^{3}T_{3}^{2}-603/
272T_{1}^{8}T_{2}^{2}T_{3}^{2}+301/136T_{1}^{6}T_{2}^{2}T_{3}^{2}-373/136T_{1}^{4}T_{2}^{
2}T_{3}^{2}+129/34T_{1}^{2}T_{2}^{2}T_{3}^{2}-47/17T_{2}^{2}T_{3}^{2}+53/68T_{1}^{9}T_{2}T_{
3}^{2}-273/136T_{1}^{7}T_{2}T_{3}^{2}+125/68T_{1}^{5}T_{2}T_{3}^{2}-249/68T_{1}^{3}T_{2}T_{
3}^{2}+117/34T_{1}T_{2}T_{3}^{2}+111/136T_{1}^{8}T_{3}^{2}+209/136T_{1}^{6}T_{3}^{2}+3/34T_{
1}^{4}T_{3}^{2}-263/68T_{1}^{2}T_{3}^{2}+60/17T_{3}^{2}+235/136T_{1}^{8}T_{2}^{3}T_{3}-147/
68T_{1}^{6}T_{2}^{3}T_{3}+91/34T_{1}^{4}T_{2}^{3}T_{3}-67/34T_{1}^{2}T_{2}^{3}T_{3}+37/34T_{
2}^{3}T_{3}-297/136T_{1}^{9}T_{2}^{2}T_{3}+199/68T_{1}^{7}T_{2}^{2}T_{3}-203/68T_{1}^{5}T_{
2}^{2}T_{3}+95/17T_{1}^{3}T_{2}^{2}T_{3}-94/17T_{1}T_{2}^{2}T_{3}-115/34T_{1}^{8}T_{2}T_{3}+
1/34T_{1}^{6}T_{2}T_{3}-83/34T_{1}^{4}T_{2}T_{3}+173/17T_{1}^{2}T_{2}T_{3}-127/17T_{2}T_{3}
-1/34T_{1}^{9}T_{3}+15/68T_{1}^{7}T_{3}+8/17T_{1}^{5}T_{3}-29/34T_{1}^{3}T_{3}-5/17T_{1}T_{
3}+1/4T_{1}^{9}T_{2}^{3}-89/136T_{1}^{7}T_{2}^{3}+95/68T_{1}^{5}T_{2}^{3}-249/68T_{1}^{3}T_{
2}^{3}+123/34T_{1}T_{2}^{3}+145/68T_{1}^{8}T_{2}^{2}+59/136T_{1}^{6}T_{2}^{2}+11/4T_{1}^{
4}T_{2}^{2}-455/68T_{1}^{2}T_{2}^{2}+70/17T_{2}^{2}-1/2T_{1}^{9}T_{2}-5/68T_{1}^{7}T_{2}-111/
68T_{1}^{5}T_{2}+67/34T_{1}^{3}T_{2}-7/34T_{1}T_{2}+44/17T_{1}^{8}+165/68T_{1}^{6}-109/68T_{
1}^{4}-63/17T_{1}^{2}+171/34$, \\
$\widehat{f_8} = T_{8}+55/136T_{1}^{8}T_{2}^{3}T_{3}^{2}+5/17T_{1}^{6}T_{2}^{3}T_{3}^{2}+9/8T_{1}^{4}T_{2}^{
3}T_{3}^{2}-227/68T_{1}^{2}T_{2}^{3}T_{3}^{2}+183/68T_{2}^{3}T_{3}^{2}-27/34T_{1}^{9}T_{2}^{
2}T_{3}^{2}+65/68T_{1}^{7}T_{2}^{2}T_{3}^{2}-197/136T_{1}^{5}T_{2}^{2}T_{3}^{2}+169/68T_{1}^{
3}T_{2}^{2}T_{3}^{2}-113/68T_{1}T_{2}^{2}T_{3}^{2}-97/136T_{1}^{8}T_{2}T_{3}^{2}+111/136T_{
1}^{6}T_{2}T_{3}^{2}-42/17T_{1}^{4}T_{2}T_{3}^{2}+321/68T_{1}^{2}T_{2}T_{3}^{2}-47/17T_{2}T_{
3}^{2}+31/34T_{1}^{9}T_{3}^{2}-31/136T_{1}^{7}T_{3}^{2}+41/34T_{1}^{5}T_{3}^{2}-9/4T_{1}^{
3}T_{3}^{2}+41/34T_{1}T_{3}^{2}+3/136T_{1}^{9}T_{2}^{3}T_{3}-113/68T_{1}^{7}T_{2}^{3}T_{3}+
137/68T_{1}^{5}T_{2}^{3}T_{3}-8/17T_{1}^{3}T_{2}^{3}T_{3}-2/17T_{1}T_{2}^{3}T_{3}+219/136T_{
1}^{8}T_{2}^{2}T_{3}+31/68T_{1}^{6}T_{2}^{2}T_{3}+163/68T_{1}^{4}T_{2}^{2}T_{3}-161/17T_{1}^{
2}T_{2}^{2}T_{3}+131/17T_{2}^{2}T_{3}-1/136T_{1}^{9}T_{2}T_{3}+215/68T_{1}^{7}T_{2}T_{3}-117/
34T_{1}^{5}T_{2}T_{3}+79/34T_{1}^{3}T_{2}T_{3}-31/34T_{1}T_{2}T_{3}-199/136T_{1}^{8}T_{3}-7/
68T_{1}^{6}T_{3}-57/34T_{1}^{4}T_{3}+219/34T_{1}^{2}T_{3}-167/34T_{3}+5/68T_{1}^{8}T_{2}^{3}+
33/136T_{1}^{6}T_{2}^{3}-29/68T_{1}^{4}T_{2}^{3}+1/68T_{1}^{2}T_{2}^{3}+10/17T_{2}^{3}-167/
136T_{1}^{9}T_{2}^{2}-201/136T_{1}^{7}T_{2}^{2}+111/34T_{1}^{5}T_{2}^{2}-71/68T_{1}^{3}T_{
2}^{2}-11/17T_{1}T_{2}^{2}+9/8T_{1}^{8}T_{2}+41/34T_{1}^{6}T_{2}-39/68T_{1}^{4}T_{2}-101/
34T_{1}^{2}T_{2}+59/17T_{2}+89/68T_{1}^{9}+127/68T_{1}^{7}-217/68T_{1}^{5}+15/17T_{1}^{3}+19/
34T_{1}$, \\
$\widehat{f_9} = T_{9}-253/544T_{1}^{8}T_{2}^{3}T_{3}^{2}+19/16T_{1}^{6}T_{2}^{3}T_{3}^{2}+71/136T_{1}^{4}T_{
2}^{3}T_{3}^{2}-409/136T_{1}^{2}T_{2}^{3}T_{3}^{2}+245/136T_{2}^{3}T_{3}^{2}+163/272T_{1}^{
9}T_{2}^{2}T_{3}^{2}+193/272T_{1}^{7}T_{2}^{2}T_{3}^{2}-121/136T_{1}^{5}T_{2}^{2}T_{3}^{2}
-53/136T_{1}^{3}T_{2}^{2}T_{3}^{2}+20/17T_{1}T_{2}^{2}T_{3}^{2}+57/272T_{1}^{8}T_{2}T_{3}^{2}
-183/136T_{1}^{6}T_{2}T_{3}^{2}-121/136T_{1}^{4}T_{2}T_{3}^{2}+59/68T_{1}^{2}T_{2}T_{3}^{2}
-18/17T_{2}T_{3}^{2}-105/272T_{1}^{9}T_{3}^{2}-5/34T_{1}^{7}T_{3}^{2}+23/136T_{1}^{5}T_{3}^{
2}+37/34T_{1}^{3}T_{3}^{2}+1/17T_{1}T_{3}^{2}-209/272T_{1}^{9}T_{2}^{3}T_{3}-4/17T_{1}^{7}T_{
2}^{3}T_{3}+213/136T_{1}^{5}T_{2}^{3}T_{3}+19/17T_{1}^{3}T_{2}^{3}T_{3}-33/17T_{1}T_{2}^{
3}T_{3}+61/68T_{1}^{8}T_{2}^{2}T_{3}+183/136T_{1}^{6}T_{2}^{2}T_{3}-287/68T_{1}^{4}T_{2}^{
2}T_{3}+71/68T_{1}^{2}T_{2}^{2}T_{3}+37/17T_{2}^{2}T_{3}+237/136T_{1}^{9}T_{2}T_{3}+157/
136T_{1}^{7}T_{2}T_{3}-15/34T_{1}^{5}T_{2}T_{3}-209/68T_{1}^{3}T_{2}T_{3}+35/34T_{1}T_{2}T_{
3}-29/34T_{1}^{8}T_{3}-81/34T_{1}^{6}T_{3}+207/68T_{1}^{4}T_{3}-39/34T_{1}^{2}T_{3}+3/34T_{
3}+4/17T_{1}^{8}T_{2}^{3}-1/2T_{1}^{6}T_{2}^{3}+311/136T_{1}^{4}T_{2}^{3}-25/17T_{1}^{2}T_{
2}^{3}+3/68T_{2}^{3}-261/272T_{1}^{9}T_{2}^{2}+107/136T_{1}^{7}T_{2}^{2}-179/136T_{1}^{5}T_{
2}^{2}+21/68T_{1}^{3}T_{2}^{2}+5/17T_{1}T_{2}^{2}-7/272T_{1}^{8}T_{2}+307/136T_{1}^{6}T_{2}
-83/34T_{1}^{4}T_{2}+257/68T_{1}^{2}T_{2}-171/68T_{2}+57/136T_{1}^{9}-209/136T_{1}^{7}-3/
17T_{1}^{5}-7/68T_{1}^{3}-25/34T_{1}$, \\
$\widehat{f_{10}} = T_{10}+15/136T_{1}^{8}T_{2}^{3}T_{3}^{2}+29/68T_{1}^{6}T_{2}^{3}T_{3}^{2}-3/136T_{1}^{4}T_{
2}^{3}T_{3}^{2}-24/17T_{1}^{2}T_{2}^{3}T_{3}^{2}+43/34T_{2}^{3}T_{3}^{2}+27/68T_{1}^{9}T_{
2}^{2}T_{3}^{2}+43/272T_{1}^{7}T_{2}^{2}T_{3}^{2}-93/68T_{1}^{5}T_{2}^{2}T_{3}^{2}+203/136T_{
1}^{3}T_{2}^{2}T_{3}^{2}-3/17T_{1}T_{2}^{2}T_{3}^{2}-9/34T_{1}^{8}T_{2}T_{3}^{2}-25/68T_{1}^{
6}T_{2}T_{3}^{2}-47/68T_{1}^{4}T_{2}T_{3}^{2}+177/68T_{1}^{2}T_{2}T_{3}^{2}-61/34T_{2}T_{3}^{
2}-39/136T_{1}^{9}T_{3}^{2}-29/136T_{1}^{7}T_{3}^{2}+207/136T_{1}^{5}T_{3}^{2}-5/4T_{1}^{
3}T_{3}^{2}+11/68T_{1}T_{3}^{2}+23/272T_{1}^{9}T_{2}^{3}T_{3}-11/68T_{1}^{7}T_{2}^{3}T_{3}+
19/136T_{1}^{5}T_{2}^{3}T_{3}+1/34T_{1}^{3}T_{2}^{3}T_{3}-2/17T_{1}T_{2}^{3}T_{3}+33/68T_{
1}^{8}T_{2}^{2}T_{3}-23/136T_{1}^{6}T_{2}^{2}T_{3}+5/34T_{1}^{4}T_{2}^{2}T_{3}-83/68T_{1}^{
2}T_{2}^{2}T_{3}+41/34T_{2}^{2}T_{3}-9/68T_{1}^{9}T_{2}T_{3}+107/136T_{1}^{7}T_{2}T_{3}-15/
34T_{1}^{5}T_{2}T_{3}+39/68T_{1}^{3}T_{2}T_{3}-7/17T_{1}T_{2}T_{3}-57/68T_{1}^{8}T_{3}-41/
68T_{1}^{6}T_{3}+5/68T_{1}^{4}T_{3}+33/17T_{1}^{2}T_{3}-65/34T_{3}-45/272T_{1}^{8}T_{2}^{3}+
3/136T_{1}^{6}T_{2}^{3}-69/68T_{1}^{4}T_{2}^{3}+67/34T_{1}^{2}T_{2}^{3}-59/68T_{2}^{3}-29/
136T_{1}^{9}T_{2}^{2}-27/34T_{1}^{7}T_{2}^{2}+267/136T_{1}^{5}T_{2}^{2}-30/17T_{1}^{3}T_{2}^{
2}+37/68T_{1}T_{2}^{2}+27/34T_{1}^{8}T_{2}+15/68T_{1}^{6}T_{2}+105/68T_{1}^{4}T_{2}-117/34T_{
1}^{2}T_{2}+40/17T_{2}+9/34T_{1}^{9}+109/136T_{1}^{7}-97/68T_{1}^{5}+77/68T_{1}^{3}-21/34T_{
1}$.
\end{example}

\begin{example}\label{gal5}
$G = 11T1 = C(11) = 11$, 
$p=47$, \\$x=(4, 10, 23, 5, 30, 19, 16, 26, 34, 6, 14)$, $d(f)=41426511213649$, \\
$\widehat{f_1} = T_{1}^{11}+T_{1}^{10}-10T_{1}^{9}-9T_{1}^{8}+36T_{1}^{7}+28T_{1}^{6}-56T_{1}^{5}-35T_{1}^{4}+
35T_{1}^{3}+15T_{1}^{2}-6T_{1}-1$, \\
$\widehat{f_2} = T_{2}+2T_{1}^{6}+3T_{1}^{5}-8T_{1}^{4}-11T_{1}^{3}+6T_{1}^{2}+7T_{1}$, \\
$\widehat{f_3} = T_{3}+T_{1}^{8}-7T_{1}^{6}+15T_{1}^{4}-T_{1}^{3}-11T_{1}^{2}+T_{1}+2$, \\
$\widehat{f_4} = T_{4}+T_{1}^{9}+T_{1}^{8}-8T_{1}^{7}-8T_{1}^{6}+19T_{1}^{5}+19T_{1}^{4}-12T_{1}^{3}-12T_{1}^{2}$, \\
$\widehat{f_5} = T_{5}-T_{1}^{9}-2T_{1}^{8}+5T_{1}^{7}+11T_{1}^{6}-7T_{1}^{5}-19T_{1}^{4}+3T_{1}^{3}+12T_{1}^{2}
-2$, \\
$\widehat{f_6} = T_{6}-T_{1}^{10}-2T_{1}^{9}+6T_{1}^{8}+14T_{1}^{7}-8T_{1}^{6}-30T_{1}^{5}-4T_{1}^{4}+21T_{1}^{
3}+7T_{1}^{2}-4T_{1}-1$, \\
$\widehat{f_7} = T_{7}+T_{1}^{9}+2T_{1}^{8}-4T_{1}^{7}-10T_{1}^{6}+T_{1}^{5}+14T_{1}^{4}+7T_{1}^{3}-5T_{1}^{2}
-2T_{1}+1$, \\
$\widehat{f_8} = T_{8}+T_{1}^{10}+2T_{1}^{9}-6T_{1}^{8}-13T_{1}^{7}+10T_{1}^{6}+26T_{1}^{5}-4T_{1}^{4}-18T_{1}^{
3}-T_{1}^{2}+3T_{1}+1$, \\
$\widehat{f_9} = T_{9}-T_{1}^{8}-2T_{1}^{7}+3T_{1}^{6}+9T_{1}^{5}+2T_{1}^{4}-11T_{1}^{3}-8T_{1}^{2}+2T_{1}+1$, \\
$\widehat{f_{10}} = T_{10}-T_{1}^{8}-T_{1}^{7}+6T_{1}^{6}+6T_{1}^{5}-9T_{1}^{4}-8T_{1}^{3}+3T_{1}^{2}+2T_{1}$, \\
$\widehat{f_{11}} = T_{11}-T_{1}^{9}+9T_{1}^{7}+T_{1}^{6}-27T_{1}^{5}-6T_{1}^{4}+30T_{1}^{3}+9T_{1}^{2}-9T_{1}-2$.
\end{example}

\begin{example}\label{gal6}
$G = 13T1 = D(13)=13:2$, 
$p=107$, \\$x=(105, 43, 77, 92, 30, 10, 78, 95, 12, 65, 41, 44, 58)$, 
\\$d(f)=49198772714708123987759220407480370279361$, \\
$\widehat{f_1} = T_{1}^{13}-T_{1}^{12}-24T_{1}^{11}+19T_{1}^{10}+190T_{1}^{9}-116T_{1}^{8}-601T_{1}^{7}+246T_{
1}^{6}+738T_{1}^{5}-215T_{1}^{4}-291T_{1}^{3}+68T_{1}^{2}+10T_{1}-1$, \\
$\widehat{f_2} = T_{2}+39245710/435105007T_{1}^{12}-4017330/435105007T_{1}^{11}-663999688/435105007T_{1}^{10}
-15607490/435105007T_{1}^{9}+3456904864/435105007T_{1}^{8}+134667006/435105007T_{1}^{7}
-5698557225/435105007T_{1}^{6}+447188240/435105007T_{1}^{5}+1889740853/435105007T_{1}^{4}
-270079487/435105007T_{1}^{3}+727404837/435105007T_{1}^{2}-157270291/435105007T_{1}-28180094/
435105007$, \\
$\widehat{f_3} = T_{3}-4746020/435105007T_{1}^{12}-6313682/435105007T_{1}^{11}-6065045/435105007T_{1}^{10}+
57287764/435105007T_{1}^{9}+603646610/435105007T_{1}^{8}-233263773/435105007T_{1}^{7} \\
-2181047898/435105007T_{1}^{6}-610031743/435105007T_{1}^{5}+1754456916/435105007T_{1}^{4}+
720178584/435105007T_{1}^{3}-537150984/435105007T_{1}^{2}-93115116/435105007T_{1} \\
+46897278/435105007$, \\
$\widehat{f_4} = T_{4}+9822090/435105007T_{1}^{12}-14190756/435105007T_{1}^{11}-301466524/435105007T_{1}^{10}+
199349313/435105007T_{1}^{9}+6553388/1372571T_{1}^{8}-1033567056/435105007T_{1}^{7} \\
-4976432274/435105007T_{1}^{6}+2872495821/435105007T_{1}^{5}+4352154340/435105007T_{1}^{4}-2668274931/
435105007T_{1}^{3}-852807980/435105007T_{1}^{2}+140070067/435105007T_{1}+9106046/435105007$, \\
$\widehat{f_5} = T_{5}-1658580/435105007T_{1}^{12}-26299035/435105007T_{1}^{11}+12325374/435105007T_{1}^{10}+
224244824/435105007T_{1}^{9}-204700811/435105007T_{1}^{8}-324659371/435105007T_{1}^{7}+
1111594779/435105007T_{1}^{6}+169516123/435105007T_{1}^{5}-1996636088/435105007T_{1}^{4}+
243678450/435105007T_{1}^{3}+10746022/5242229T_{1}^{2}-252055005/435105007T_{1} \\
+13891072/435105007$, \\
$\widehat{f_6} = T_{6}-403010/18917609T_{1}^{12}+49410411/435105007T_{1}^{11}+61393444/435105007T_{1}^{10}
-411495261/435105007T_{1}^{9}-179743148/435105007T_{1}^{8}+1332664360/435105007T_{1}^{7}+
453161115/435105007T_{1}^{6}-2395377948/435105007T_{1}^{5}-1354841892/435105007T_{1}^{4}+
1052235225/435105007T_{1}^{3}+701947707/435105007T_{1}^{2}-7025266/435105007T_{1} \\
-12442298/435105007$, \\
$\widehat{f_7} = T_{7}-12088160/435105007T_{1}^{12}+36013925/435105007T_{1}^{11}+240777437/435105007T_{1}^{10}
-687864750/435105007T_{1}^{9}-1337098960/435105007T_{1}^{8}+3736584070/435105007T_{1}^{7}+
1576974451/435105007T_{1}^{6}-5658657766/435105007T_{1}^{5}+836138773/435105007T_{1}^{4}+
2352823294/435105007T_{1}^{3}-30848390/18917609T_{1}^{2}-130536919/435105007T_{1} \\
+4661965/435105007$, \\
$\widehat{f_8} = T_{8}-19622855/435105007T_{1}^{12}-42521665/435105007T_{1}^{11}+348666670/435105007T_{1}^{10}+
708283196/435105007T_{1}^{9}-1857480891/435105007T_{1}^{8}-3548072906/435105007T_{1}^{7}+
3107626693/435105007T_{1}^{6}+5709162399/435105007T_{1}^{5}-1376716168/435105007T_{1}^{4}
-2274067840/435105007T_{1}^{3}+241517773/435105007T_{1}^{2}+157457552/435105007T_{1}-33657606/
435105007$, \\
$\widehat{f_9} = T_{9}-15794040/435105007T_{1}^{12}+21092920/435105007T_{1}^{11}+308767305/435105007T_{1}^{10}
-291624628/435105007T_{1}^{9}-1982697657/435105007T_{1}^{8}+1898428538/435105007T_{1}^{7}+
3953390793/435105007T_{1}^{6}-4188193592/435105007T_{1}^{5}-1286268704/435105007T_{1}^{4}+
2373302361/435105007T_{1}^{3}-681114322/435105007T_{1}^{2}+5052302/435105007T_{1} \\
+13080401/435105007$, \\
$\widehat{f_{10}} = T_{10}-199965/435105007T_{1}^{12}-18647220/435105007T_{1}^{11}+101510623/435105007T_{1}^{10}+
213783862/435105007T_{1}^{9}-959373536/435105007T_{1}^{8}-596249442/435105007T_{1}^{7}+
2642799996/435105007T_{1}^{6}+322557418/435105007T_{1}^{5}-3188331064/435105007T_{1}^{4}+
183121060/435105007T_{1}^{3}+1334545000/435105007T_{1}^{2}-223977170/435105007T_{1}-11326118/
435105007$, \\
$\widehat{f_{11}} = T_{11}-1504028/435105007T_{1}^{12}-6306587/435105007T_{1}^{11}+33653768/435105007T_{1}^{10}+
49785758/435105007T_{1}^{9}-4491350/435105007T_{1}^{8}-71935693/435105007T_{1}^{7}\\
-32802386/18917609T_{1}^{6}-1263306769/435105007T_{1}^{5}+1818495180/435105007T_{1}^{4}+1726453403/
435105007T_{1}^{3}-1233104652/435105007T_{1}^{2}+254687554/435105007T_{1}-13595364/435105007$, \\
$\widehat{f_{12}} = T_{12}+28679195/435105007T_{1}^{12}-15814791/435105007T_{1}^{11}-295797896/435105007T_{1}^{10}+
201896280/435105007T_{1}^{9}+596276584/435105007T_{1}^{8}-742446913/435105007T_{1}^{7}+
755441993/435105007T_{1}^{6}+1727640755/435105007T_{1}^{5}-1464470476/435105007T_{1}^{4}
-1242319218/435105007T_{1}^{3}+416242961/435105007T_{1}^{2}+145628554/435105007T_{1}-12605435/
435105007$, \\
$\widehat{f_{13}} = T_{13}-12864117/435105007T_{1}^{12}+27593810/435105007T_{1}^{11}+160234532/435105007T_{1}^{10}
-248038868/435105007T_{1}^{9}-208665701/435105007T_{1}^{8}-552148820/435105007T_{1}^{7}+9502455/
435105007T_{1}^{6}+2867007062/435105007T_{1}^{5}+16278330/435105007T_{1}^{4} \\
-2197050901/435105007T_{1}^{3}-299887196/435105007T_{1}^{2}+161083738/435105007T_{1}+24170153/435105007$.
\end{example}

Note that examples \ref{gal1},\ref{gal2}, \ref{gal3}, \ref{gal4}, \ref{gal5} and \ref{gal6} serve as examples of 
our constructive version of the theorem of Galois, since the respective groups are simple. 

Now let us say a word about the running time of our algorithm for computing of the polynomials $\widehat{f_i}$. 
The following table compares the running time (in milliseconds) of KANT that is needed for computing all the preliminaries
(the prime $p$ such that $f$ splits into disjoint linear factors over $\Q_p$, 
the Galois group of $f$, and approximations of the $p$-adic zeros of $f$ with exactly the precision we need, 
such that the Galois group permutes the zeros) on the one hand and the running time (in milliseconds) 
of the GAP programme which from these results computes the polynomials $\widehat{f_i}$, 
by \eqref{lag}, on the other hand. The second line of the table shows the degree of the respective Galois group. 
Clearly the running times of the algorithms are roughly proportional to those degrees. 
All computations were done on a machine with 2 Intel processors, working at 2 GHz each, and 2 GB memory, 
running under Linux. 

\begin{center}
\begin{tabular}{|c|c|c|c|c|c|c|c|c|c} \hline 
group  & 5T1 & 5T3 & 6T3 & 6T5 & 7T1 & 7T3 & 8T20 & 8T24 & 8T32  \\ \hline 
degree & 5   & 20  & 12  & 18  & 7   & 21  & 32   & 48   & 96    \\ \hline 
KANT   & 90  & 140 & 360 & 160 & 180 & 280 & 450  & 580  & 900   \\ \hline 
GAP    & 30  & 140 & 90  & 160 & 450 & 810 & 2860 & 2130 & 12580 \\ \hline 
\end{tabular}
\end{center}

\begin{center}
\begin{tabular}{c|c|c|c|c|c|c|c|c|c|} \hline 
\, & 9T7   & 9T11  & 9T18  & 10T5  & 10T8  & 10T11 & 11T1  & 13T1  \\ \hline 
\, & 27    & 54    & 108   & 40    & 80    & 120   & 11    & 13    \\ \hline 
\, & 880   & 1690  & 1520  & 950   & 1050  & 6200  &       &       \\ \hline 
\, & 18900 & 12310 & 29250 & 15060 & 90480 & 71220 & 44570 & 39890 \\ \hline 
\end{tabular}
\end{center}

The relatively long running time of the KANT computation is mainly due to search for the prime $p$ that we want. 
Compared to that, the Galois group is computed relatively quickly. Of course, one could also embed $L$ into 
another $p$-adic field and thereby get rid of the search for the paricular $p$ that we have been working with. 
There is another reason why it would be good to work with a different $p$ than the one we have been working with:
As J\"urgen Kl\"uners told me, for primes $p$ as we chose them (i.e. such that $f$ splits into disjoint linear factors), 
the KANT algorithm for computing the Galois group will tend to infinity as the size of the group gets bigger, 
see also Section 4 in \cite{klue}. The reason why I preferred to work with my $p$ is to keep the 
technicalities of Secion \ref{secnum} as easy as possible. The disadvantage of doing so is that
for the polynomials of degree $\geq11$, the KANT algorithm I used for computing the Galois group 
does not yield the group with absolute certainty, but only with very high probability. 
This is why I decided not to write down the running time of the KANT algorithm in these cases. 

The selection of examples may seem quite limited since the degrees of the field extensions are not very high. 
However, I tested the algorithms for a much wider range of examples. The next table presents the running times 
of the algorithms, for some sample groups. For the polynomials having the respective groups, 
I again consulted the Kl\"uners--Malle database. For each group they give a selection of appropriate polynomials, 
of which I always took the first. 

\begin{center}
\begin{tabular}{|c|c|c|c|c|c|c|c|c} \hline 
group  & 10T8  & 8T32  & 7T5  & 8T42  & 9T25  & 6T15  & 10T28  & 9T27   \\ \hline 
degree & 80    & 96    & 168  & 288   & 324   & 360   & 400    & 504    \\ \hline 
KANT   & 1050  & 900   & 1520 & 970   & 1490  & 2820  & 7300   & 15000  \\ \hline 
GAP    & 90480 & 12580 & 9730 & 20080 & 67930 & 14580 & 909330 & 397340 \\ \hline 
\end{tabular}
\end{center}

\begin{center}
\begin{tabular}{c|c|c|c|c|c|c|c|c|} \hline 
\, & 10T33   & 8T47   & 9T31   & 9T32    & 10T36   & 7T6     & 10T39    \\ \hline 
\, & 800     & 1152   & 1296   & 1512    & 1920    & 2520    & 3840     \\ \hline 
\, & 17000   & 24000  & 140000 & 16000   & 23      & 110000  & 23000    \\ \hline 
\, & 2402790 & 238340 & 403390 & 1967120 & 7366290 & 1508310 & 22466640 \\ \hline 
\end{tabular}
\end{center}

\section{Questions}\label{secque}
In Section \ref{secnum}, at a certain point we multiplied $\widehat{f_i}$ with the factor $\Delta_i$ (essentially 
a power of the discriminant) in order to obtain a polynomial with coefficients in $\Z$. Therefore, as for the rational 
polynomial $\widehat{f_i}$, one would expect that the denominators that occur in its coefficients are
in the magnitude of $\Delta_i$. But in all examples computed so far, the denominators are significantly smaller
than $\Delta_i$. This phenomenon can be explained by a very elementary argument in the case 
where $f$ has the Galois group $S_n$, $n$ being the degree of $f$. 
In order to do so, we recall the definition of $K_i$ and $f_i$. 

We have $K_i=K(x_1,\ldots,x_i)$ for all $i$, and $f_i$ is the minimal polynomial of $x_i$ over $K_{i-1}$. 
Clearly $x_i$ is a zero of $f$, hence $f_i$ appears as a factor of $f$ when factorising this polynomial 
over $K_{i-1}$. We can even characterise $f_i$ as being the monic irreducible factor 
of $f(T_i)\in K_{i-1}[T_i]$ having $x_i$ as a zero. (As already mentioned in Section \ref{secgen}, 
this is the way $f_i$ is characterised in \cite{anai}, and this is the basis upon \cite{anai} is built up.)

Now let us study the case $G=S_n$. There the degree of the field extension $L|\Q$ equals $n!$. 
The degree $d_i$ of the field extension $K_i|K_{i-1}$ equals the degree of $f_i$, a factor of $f$. 
Over $K_{i-1}$, the polynomial $f$ has the factors $T-x_j$, $j\leq i-1$. Hence the degree of 
the factor $f_i$ is bounded above by $n-i+1$. From $\prod_{j\leq n}d_j$, we conclude that $d_i=n-i+1$. 
Therefore the polynomial $f$ has a factorisation $f(T_i)=(T_i-x_1)\ldots(T_i-x_{i-1})f_i(T_i)$ over $K_{i-1}[T_i]$. 
From this one can easily derive an explicit formula for $f_i(T_i)$ in terms of the coefficients of $f$ 
and of $x_1,\ldots,x_{i-1}$. It turns out that $f_i(T_i)$ is a polynomial expression in the coefficients of $f$, 
in $x_1,\ldots,x_{i-1}$ and in $T_i$. Thus $\widehat{f_i}$ is a polynomial expression 
in the coefficients of $f$ and in $T_1,\ldots,T_i$. If in particular the coefficients of $f$ are integers, 
it follows that $\widehat{f_i}$ is an integer polynomial in $T_1,\ldots,T_i$. 

However, in the case when $G$ is a proper subset of $S_n$, the question for a fair upper bound 
for the denominators of $\widehat{f_i}$ seems to be rather difficult. 

Let us ask some more questions, concerning a generalisation of Theorem \ref{buch}. 
The Buchberger--M\"oller algorithm serves to compute the ideal $I(X)$ 
for a given set of points $X$ in $\A^n(F)$. In Theorem \ref{buch} we presented a situation 
in which the Buchberger--M\"oller algorithm is actually obsolete since 
one can compute the Gr\"obner basis of $I(X)$ by using the explicit formula \eqref{laggen} -- 
at least for lexicographic ordering. Now the obvious question is: 
Is there an analogous formula also for other monomial orderings (as classified in \cite{rob85} and \cite{rob86})? 
Say, as a starting point, we could replace the set $X=\{(\sigma(x_1),\ldots,\sigma(x_n));\sigma\in G\}$ by some
$X=\{(\sigma(x_1),\ldots,\sigma(x_n))A;\sigma\in G\}$ for an appropriate matrix of size $n$ oder $K$. 
(Example: $A_{i,j}=\delta_{n-i,j}$ will lead us to the lexicographical ordering, 
with the indeterminates ordered the other way round.)
What about a different kind of action of $G$ on $X$? 
(Example: The one-dimensional Lagrange interpolation, 
in which there is no $G$ acting as in Theorem \ref{buch},
yet there definitely is an analogue of \eqref{laggen}.)

\section*{Acknowledgements}

I wish to thank Herwig Hauser for awakening my interest in the relation ideal and Kurt Girstmair for telling 
me what this ideal is good for and for his help during my work on this paper. 
Furthermore, I wish to thank Alexander Ostermann for his help on interpolation and for the great motivation, 
J\"{u}rgen Kl\"{u}ners for his help in KANT questions, and Leonhard Wieser for introducing me 
to the Buchberger--M\"oller algorithm. It was from his excellent diploma thesis \cite{leo} 
that I learned about this algorithm. Special thanks to Michael Redinger for his constant support 
in software problems of any kind!
Finally, I wish to thank Fran\c{c}ois Brunault for many very useful comments on the text, 
and my brother Thomas Lederer for helping to improve my English.

\bibliographystyle{elsart-harv}
\bibliography{relation_buch}

\begin{thebibliography}{21}
\expandafter\ifx\csname natexlab\endcsname\relax\def\natexlab#1{#1}\fi
\expandafter\ifx\csname url\endcsname\relax
  \def\url#1{\texttt{#1}}\fi
\expandafter\ifx\csname urlprefix\endcsname\relax\def\urlprefix{URL }\fi

\bibitem[{Abbott et~al.(2000)Abbott, Bigatti, Kreuzer, and Robbiano}]{abbott}
Abbott, J., Bigatti, A., Kreuzer, M., Robbiano, L., 2000. Computing ideals of
  points. J. Symbolic Comput. 30~(4), 341--356.

\bibitem[{Alonso et~al.(2003)Alonso, Marinari, and Mora}]{alonso}
Alonso, M.~E., Marinari, M.~G., Mora, T., 2003. The big mother of all
  dualities: {M}\"oller algorithm. Comm. Algebra 31~(2), 783--818.

\bibitem[{Anai et~al.(1996)Anai, Noro, and Yokoyama}]{anai}
Anai, H., Noro, M., Yokoyama, K., 1996. Computation of the splitting fields and
  the {G}alois groups of polynomials. In: Algorithms in algebraic geometry and
  applications (Santander, 1994). Vol. 143 of Progr. Math. Birkh\"auser, Basel,
  pp. 29--50.

\bibitem[{Aubry and Valibouze(2000)}]{aub}
Aubry, P., Valibouze, A., 2000. Using {G}alois ideals for computing relative
  resolvents. J. Symbolic Comput. 30~(6), 635--651, algorithmic methods in
  Galois theory.

\bibitem[{Becker and Weispfenning(1993)}]{beck}
Becker, T., Weispfenning, V., 1993. Gr\"obner bases. Vol. 141 of Graduate Texts
  in Mathematics. Springer-Verlag, New York, a computational approach to
  commutative algebra, In cooperation with Heinz Kredel.

\bibitem[{Buchberger(1965)}]{buch}
Buchberger, B., 1965. Ein {A}lgorithmus zum {A}uffinden der {B}asiselemente des
  {R}estklassenrings nach einem nulldimensionalen {P}olynomideal. Ph.D. thesis,
  Universit{\"{a}}t Innsbruck.

\bibitem[{Cohn(1989)}]{cohn}
Cohn, P.~M., 1989. Algebra. {V}ol. 2. John Wiley \& Sons Ltd., Chichester.

\bibitem[{Daberkow et~al.(1997)Daberkow, Fieker, Kl{\"u}ners, Pohst, Roegner,
  Sch{\"o}rnig, and Wildanger}]{KANT}
Daberkow, M., Fieker, C., Kl{\"u}ners, J., Pohst, M., Roegner, K.,
  Sch{\"o}rnig, M., Wildanger, K., 1997. K{ANT} {V}4. J. Symbolic Comput.
  24~(3-4), 267--283.

\bibitem[{GAP(2004)}]{GAP4}
GAP, 2004. {GAP -- Groups, Algorithms, and Programming, Version 4.4}. The
  GAP~Group, \verb+(http://www.gap-system.org)+.

\bibitem[{Geissler and Kl{\"u}ners(2000)}]{klue}
Geissler, K., Kl{\"u}ners, J., 2000. Galois group computation for rational
  polynomials. J. Symbolic Comput. 30~(6), 653--674, algorithmic methods in
  Galois theory.

\bibitem[{Gr{\"o}bner(1970)}]{gro}
Gr{\"o}bner, W., 1970. Algebraische {G}eometrie. {$2$}. {T}eil: {A}rithmetische
  {T}heorie der {P}olynomringe. Bibliographisches Institut, Mannheim, b. I.
  Hochschultaschenb\"ucher, 737/737a$\sp{\ast}$.

\bibitem[{Huppert(1967)}]{hup}
Huppert, B., 1967. Endliche {G}ruppen. {I}. Die Grundlehren der Mathematischen
  Wissenschaften, Band 134. Springer-Verlag, Berlin.

\bibitem[{Lederer(2002)}]{led}
Lederer, M., 2002. Explizite {K}onstruktionen in
  {Z}erf{\"{a}}llungsk{\"{o}}rpern von {P}olynomen. Master's thesis,
  Universit{\"{a}}t Innsbruck.

\bibitem[{Lederer(2004)}]{parma}
Lederer, M., 2004. Explicit constructions in splitting fields of polynomials.
  Riv. Mat. Univ. Parma (7) 3*, to appear.

\bibitem[{McKay and Stauduhar(1997)}]{mckay}
McKay, J., Stauduhar, R., 1997. Finding relations among the roots of an
  irreducible polynomial. In: Proceedings of the 1997 International Symposium
  on Symbolic and Algebraic Computation (Kihei, HI). ACM, New York, pp. 75--77
  (electronic).

\bibitem[{M{\"o}ller and Buchberger(1982)}]{buchbergermoeller}
M{\"o}ller, H.~M., Buchberger, B., 1982. The construction of multivariate
  polynomials with preassigned zeros. In: Computer algebra (Marseille, 1982).
  Vol. 144 of Lecture Notes in Comput. Sci. Springer, Berlin, pp. 24--31.

\bibitem[{Neukirch(1999)}]{neu}
Neukirch, J., 1999. Algebraic number theory. Vol. 322 of Grundlehren der
  Mathematischen Wissenschaften [Fundamental Principles of Mathematical
  Sciences]. Springer-Verlag, Berlin.

\bibitem[{Robbiano(1985)}]{rob85}
Robbiano, L., 1985. Term orderings on the polynomial ring. In: EUROCAL '85,
  Vol.\ 2 (Linz, 1985). Vol. 204 of Lecture Notes in Comput. Sci. Springer,
  Berlin, pp. 513--517.

\bibitem[{Robbiano(1986)}]{rob86}
Robbiano, L., 1986. On the theory of graded structures. J. Symbolic Comput.
  2~(2), 139--170.

\bibitem[{Robbiano(1998)}]{rob}
Robbiano, L., 1998. Gr\"obner bases and statistics. In: Gr\"obner bases and
  applications (Linz, 1998). Vol. 251 of London Math. Soc. Lecture Note Ser.
  Cambridge Univ. Press, Cambridge, pp. 179--204.

\bibitem[{Wieser(2004)}]{leo}
Wieser, L., 2004. Der {B}uchberger--{M}{\"{o}}ller--{A}lgorithmus und
  polynomiale {I}nterpolation in mehreren {V}ariablen. Master's thesis,
  Universit{\"{a}}t Innsbruck.

\end{thebibliography}

\end{document}